\documentclass[11pt]{article}

\usepackage{amssymb,amsmath,euscript}

\newtheorem{proposition}{Proposition}[section]
\newtheorem{lemma}[proposition]{Lemma}
\newtheorem{theorem}[proposition]{Theorem}

\def\a{\alpha}

\def\de{\delta}

\def\ga{\gamma}
\def\Ga{\Gamma}
\def\la{\lambda}

\def\si{\sigma}
\def\e{\eta}
\def\th{\theta}
\def\o{\omega}

\def\det{{\rm det}}
\def\dim{{\rm dim}_{_{\rm {\cal H}}}}

\def\Sp{{\rm Supp}}
\def\Var{{\rm Var}}

\def\sgn{{\rm sgn}}

\def\phe{pointwise H\"older exponent }

\def\N{{\mathbb N}}
\def\N{{\mathbb N}}
\def\R{{\mathbb R}}
\def\Q{{\mathbb Q}}
\def\eps{\varepsilon}

\def\E{{\mathbb E}}
\def\P{{\mathbb P}}
\def\Z{{\mathbb Z}}

\makeatletter \@addtoreset{equation}{section} \makeatother
\newcommand {\qed}%
{%
    {}\hfill
    {}\hfill
    {$\square $}%
    \vspace {0.3cm}%
    \pagebreak [2]%
    \par
}%
\newenvironment{remark}{%
    \vspace{0.3cm} \pagebreak [2]%
    \par%
    \refstepcounter{proposition}
    \noindent%
    {\bf Remark~\theproposition\  }}{\qed}%
\begin{document}

\title {Continuous Gaussian multifractional processes with random pointwise H\"older regularity (long version)}
\author{Antoine Ayache\\ Universit\'e Lille 1}

\maketitle

\begin{abstract} Let $\{X(t)\}_{t\in\R}$ be an arbitrary centered Gaussian process whose trajectories are, with probability $1$, continuous
nowhere differentiable functions. It follows from a classical result, derived from zero-one law, that, with probability
$1$, the trajectories of $X$ have the same global H\"older
regularity over any compact interval, that is the uniform H\"older exponent does not depend on the choice of a trajectory. A similar 
phenomenon happens with their local H\"older regularity measured through the local H\"older exponent. Therefore, it seems natural to ask the 
following question: does such a phenomenon also occur with their pointwise H\"older regularity measured through the pointwise H\"older 
exponent? 

In this article, using the framework of
multifractional processes, we construct a family of counterexamples
showing that the answer to this question is not always positive.
\end{abstract}

\noindent{\small{\bf Running Title}:\ Gaussian multifractional
processes with random exponents}



\section{Introduction}
\label{sec:intro} Let $\{X(t)\}_{t\in\R}$ be an arbitrary centered Gaussian
process whose trajectories are, with probability $1$, continuous 
nowhere differentiable functions over the real line $\R$. The global H\"older regularity of 
one of them, $t\mapsto X(t,\omega)$, over a non-degenerate \footnote{By non-degenerate we mean that the compact interval $J$
is not empty nor a single point.} compact
interval $J\subset\R$, is measured through the uniform H\"older
exponent $\beta_X(J,\o)$ defined as,
\begin{equation}
\label{eq:unifhold} \beta_X(J,\o)=\sup\left\{\beta \ge
0\,:\,\,\sup_{t',t''\in
J}\frac{|X(t',\o)-X(t'',\o)|}{|t'-t''|^\beta}<\infty\right\}.
 \end{equation}
The local H\"older regularity of the trajectory $t\mapsto
X(t,\omega)$ in a neighborhood of some fixed point $s\in\R$, is
measured through two different exponents: the local H\"older exponent $\widetilde{\a}_X(s,\omega)$ and the \phe $\a_X(s,\o)$. They are defined as:
\begin{equation}
\label{eq:lochold1}
\widetilde{\a}_X(s,\omega)=\sup\left\{\beta_X\big([u,v],\o\big)\,:\,\,u,v\in\R \mbox{ and } s\in (u,v)\right\}
\end{equation}
and
\begin{equation}
\label{eq:poinhold}
\a_X(s,\o)=\sup\left\{\a\ge 0\,:\,\,\limsup_{h\rightarrow 0}\frac{|X(s+h,\o)-X(s,\o)|}{|h|^\a}=0\right\}.
\end{equation}
Notice that one always has,
$$
\widetilde{\a}_X(s,\omega)\le \a_X(s,\o).
$$
Moreover, the function $s\mapsto \widetilde{\a}_X(s,\omega)$ is always lower semicontinuous over $\R$ \cite{Seuret02}, while the function 
$s\mapsto \a_X(s,\omega)$ does not necessarily satisfy such a nice property; in fact the latter function can be the liminf of any arbitrary sequence of continuous functions 
with values in $[0,1]$ (see \cite{And,DLM,J1,AJT07,Seuret02}), therefore its behavior can be quite erratic. The notion of \phe is a fundamental concept 
in the area of the multifractal analysis of deterministic and random functions \cite{J3,JMR}. It provides a sharp estimation of the asymptotic behavior of the modulus of local continuity of the function $t\mapsto X(t,\o)$ at any fixed point $s\in\R$ (see e.g. \cite{Lif} page 214 for the defintion of this modulus of continuity); indeed, (\ref{eq:poinhold}) implies that,
$$
\a_X(s,\o)=\sup\left\{\a\ge 0\,:\,\,\limsup_{\rho\rightarrow 0_+}\left(\sup\left\{\frac{|X(t,\o)-X(s,\o)|}{\rho^\a}\,:\,\, \mbox{$t\in\R$ and $|t-s|\le\rho$}\right\}\right)=0\right\}.
$$
In order to explain the main motivation behind our article, let us state the following theorem, whose proof which is given in Subsection~\ref{subsec:expdet}, implicitly relies on zero-one law. 

\begin{theorem}
\label{theo:expdet}
One denotes by $\{X(t)\}_{t\in\R}$ an arbitrary centered Gaussian process whose trajectories are, with probability $1$, continuous 
nowhere differentiable functions over the real line $\R$. The following two results hold.
\begin{itemize}
\item[(i)] For each non-degenerate compact interval $J\subset\R$, let $b_X(J)$ be the deterministic quantity defined as,
\begin{equation}
\label{eq1:expdet}
b_X(J)=\sup\left\{ b\ge 0\,:\,\,\sup_{t',t''\in
J}\frac{\E|X(t',\o)-X(t'',\o)|^2}{|t'-t''|^{2b}}<\infty\right\}.
\end{equation}
Then, one has,
\begin{equation}
\label{eq:unif-holder} \P\big\{\beta_X(J)=b_X(J)\big\}=1,
\end{equation}
where $\P\big\{\beta_X(J)=b_X(J)\big\}$ denotes the probability that the uniform H\"older exponent
$\beta_X (J)$ be equal to $b_X(J)$.
\item[(ii)] There exists $\widetilde{\Omega}$ an event of probability $1$, non depending on $s$, such that,
the local H\"older exponent $\widetilde{\a}_X$ satisfies, 
\begin{equation}
\label{eq:stongphl} 
\widetilde{\a}_X(s,\o)=\widetilde{a}_X(s), \mbox{ for all $(s,\o)\in\R\times \widetilde{\Omega}$,}
\end{equation}
where $\widetilde{a}_X(s)$ is the deterministic quantity defined as,
\begin{equation}
\label{eq2:expdet}
\widetilde{a}_X(s)=\sup\left\{b_X\big([u,v]\big)\,:\,\,u,v\in\R \mbox{ and } s\in (u,v)\right\}.
\end{equation}
\end{itemize}
\end{theorem}
Notice that Theorem~\ref{theo:expdet} Part $(ii)$, has already been obtained in \cite{HLV} (see Corollary 3.15 in this article) under the assumption that $\widetilde{a}_X$ is  continuous. This result means that, with probability $1$, the function $s\mapsto \widetilde{\a}_X(s,\omega)$ does not depend on the choice of $\omega$, whatever the centered continuous nowhere differentiable Gaussian process $\{X(t)\}_{t\in\R}$ might be. The main goal of our article is to show that, in some cases, a different phenomenon happens for the function $s\mapsto \a_X(s,\omega)$; namely we construct multifractional Gaussian processes $\{X(t)\}_{t\in\R}$ with continuous nowhere differentiable trajectories, such that with 
a strictly positive probability the function $s\mapsto \a_X(s,\omega)$ depends on the choice of $\omega$. To this end, we draw a close connection between the values of the latter function and the zero-level set $\big\{s\in\R\,:\,\,Y(s,\omega)=0\big\}$ of a Gaussian process $\{Y(s)\}_{s\in\R}$ closely related to $X$ and very similar to it.

It is worth noticing that for each $s\in\R$, there always exists a deterministic quantity $a_X(s)\in [0,1]$ such that
\begin{equation}
\label{eq:weakphe} \P\big\{\a_X(s)=a_X(s)\big\}=1.
\end{equation}
Relation (\ref{eq:weakphe}) corresponds to Lemma 3.5 in~\cite{AT03} and Proposition 6.2 in~\cite{AJT07}, it means that the deterministic function $a_X$ is 
a modification of the stochastic process $\a_X$. In view of (\ref{eq:weakphe}), the fact that, with a strictly positive probability, the function 
$s\mapsto \a_X (s,\omega)$ depends on $\o$, implies that the deterministic function $a_X$ and the stochastic process $\a_X$ are not indistinguishable; 
not indistinguishable formally  means that: for all event $\widetilde{\Omega}$ of
probability $1$, there exists $\o_0\in\widetilde{\Omega}$ and $s_0(\o_0)\in\R$, such that,
$$
\a_X \big(s_0(\o_0),\o_0\big)\ne a_X \big(s_0(\o_0) \big).
$$
In order to show that with a strictly positive probability, the function 
$s\mapsto \a_X (s,\omega)$ depends on $\o$, we use the framework of multifractional Brownian motion (mBm). Let
us now make a few recalls concerning this Gaussian process. We
denote by $\{B(t,\th)\}_{(t,\th)\in\R\times (0,1)}$ the centered Gaussian
field defined for all $(t,\th)$ as the Wiener integral,
\begin{equation}
\label{eq:fieldB}
B(t,\th)=\int_{\R}\Big\{(t-x)_{+}^{\th-1/2}-(-x)_{+}^{\th-1/2}\Big\}\,dW(x),
\end{equation}
with the convention that for every $(u,\th)\in\R^2$,
$(u)_{+}^{\th-1/2}=u^{\th-1/2}$ if $u>0$ and $(u)_{+}^{\th-1/2}=0$
else. Let $H$ be an arbitrary fixed continuous function defined on $\R$ and with
values in the open interval $(0,1)$. The mBm of
functional parameter $H$, is the centered Gaussian process
$\{X(t)\}_{t\in\R}$ defined for all $t$ as,
\begin{equation}
\label{eq:defmBm}
X(t)=B(t,H(t))=\int_{\R}\Big\{(t-x)_{+}^{H(t)-1/2}-(-x)_{+}^{H(t)-1/2}\Big\}\,dW(x).
\end{equation}
Notice that in Theorem~\ref{theo:fieldB}, we introduce a modification $\widetilde{B}$ of the field $B$ defined as a random wavelet type series (see Subsection~\ref{subsec:wav}); let us stress that, for the sake of simplicity, in our article, $B$ is often identified with $\widetilde{B}$ and the process $X$ with its modification $\widetilde{X}$, 
defined for all $t\in\R$, as $\widetilde{X}(t)=\widetilde{B}(t,H(t))$. 


MBm is an extension of fractional Brownian motion (fBm), indeed,
when the function $H$ is a constant denoted by $h$, then $X$ reduces
to a fBm of Hurst parameter $h$; the latter Gaussian process has
been widely studied since several decades, we refer to e.g.
\cite{ST94,EM02,Adler81,kahane85,Fa1} for a presentation of its
main properties. MBm was introduced, independently in \cite{PL95} and
\cite{BJR97}, to overcome an important drawback of fBm due to the
fact that its \phe remains constant all along its trajectory. Since
several years there is an increasing interest in the study of
multifractional processes (see for instance
\cite{ayache02,AT03,AJT07,ASX,Fa2,Fa3,FaL,Herbin06,MWX08,StTaqqu06,Sur}).
It has been proved in \cite{BJR97,PL95} that when $H$ is a H\"older
function over a non-degenerate compact interval $J$ and satisfies
the condition
\begin{equation}
\label{eq:condusH} \max_{t\in J} H(t)< \beta_H (J),
\end{equation}
where $\beta_H (J)$ denotes the uniform H\"older exponent of $H$
over $J$, then for all $s\in\mathring{J}$ (note that one restricts
to $\mathring{J}$, the interior of $J$, in order to avoid the border
effect), one has
\begin{equation}
\label{eq:regusX} \P\big\{\a_X(s)=H(s)\big\}=1.
\end{equation}
Later it has been shown in \cite{AJT07}, that when the condition
(\ref{eq:condusH}) is satisfied, then $\{H(s)\}_{s\in\mathring{J}}$ and $\{\a_X(s)\}_{s\in\mathring{J}}$ are indistinguishable; namely there exists
$\widetilde{\Omega}$, an event of probability $1$, non depending on
$s$ (and also non depending on $J$), such that,
\begin{equation}
\label{eq:reguniformusX} 
\a_X(s,\o)=H(s), \mbox{ for all $(s,\o)\in(\mathring{J},\widetilde{\Omega})$.}
\end{equation}
When the condition (\ref{eq:condusH}) fails to be satisfied, under the 
assumption that for all $s\in\R$ one has $\a_H (s)\neq H(s)$, one can show that, for every
$s\in\R$,
\begin{equation}
\label{eq:failcondusHlocal}\P\big\{\a_X(s)=\min\{H(s),\a_H
(s)\}\big\}=1,
\end{equation}
where $\a_H (s)$
denotes the \phe of $H$ at $s$. Observe that (\ref{eq:failcondusHlocal}) has been derived in \cite{Herbin06,HLV}, for a definition of mBm slightly different from 
(\ref{eq:defmBm}), yet the proof also works in the latter case. 

In our article, we construct examples of Gaussian
mBm's $X$ with continuous nowhere differentiable trajectories
which satisfy the following property: there exists an event $D$ of strictly positive probability,
such that for all $\o\in D$, one has,
$$
\a_X \big(s_0(\o),\o\big)\ne\min\left\{H\big(s_0(\o)\big),\a_H \big(s_0(\o)\big)\right\},
$$
for some $s_0(\o)\in\R$. In view of
(\ref{eq:failcondusHlocal}), the latter property means that the
pointwise H\"older regularity of $X$ is random, in other words, it depends on the
choice of a trajectory of $X$. Note in passing, that there are many examples of non Gaussian processes whose regularity 
is random, as for instance, discontinuous L\'evy processes \cite{J2}, multifractional 
processes with random exponent \cite{AT03,AJT07}, self-regulating
processes \cite{Bar}, or pure jump Markov processes \cite{BFJS}.

The remaining of this article is structured in the following way. 
In Section~\ref{sec:example}, we introduce a modification $\widetilde{B}$ of the  field $B$,  
having some nice properties which are useful for the
study of the pointwise H\"older regularity of the mBm $X$. Then, denoting by ${\cal J}\subseteq (0,+\infty)$ 
an arbitrary open non-empty interval, under some
condition on its parameter $H$, we show that the pointwise H\"older regularity of $\{X(t)\}_{t\in{\cal J}}$,
is closely connected with the zeros of the process
$\{Y(s)\}_{s\in {\cal J}}=\big\{(\partial_\th B)(s,H(s))\big\}_{s\in {\cal J}}$; thus it turns
out that this regularity is random, when the level set $\big\{s\in
{\cal J}\,:\,\, Y(s)=0\big\}$ is non-empty with a (strictly)
positive probability. In Section~\ref{sec:levelset}, we
prove that this is indeed the case, namely with a strictly positive probability the latter level set, is rather
large: it has a Hausdorff dimension bigger than $1-\eta-\inf_{s\in {\cal J}} H(s)$, where $\eta$ 
is a fixed strictly positive and arbitrarily small real number. Finally, some technical proofs, mainly related to wavelet methods, are given 
in Section~\ref{sec:appendice} (the Appendix).

\section{Construction of the counterexamples}
\label{sec:example} 
In order to construct continuous Gaussian multifractional Brownian motions with random pointwise H\"older regularity, first we need to show that the Gaussian random
field $\{B(t,\th)\}_{(t,\th)\in\R\times (0,1)}$ defined in (\ref{eq:fieldB}), has a modification $\{\widetilde{B}(t,\th)\}_{(t,\th)\in\R\times (0,1)}$ satisfying some 
nice properties. Namely, we need the following theorem.
\begin{theorem}
\label{theo:fieldB} 
Let $\{B(t,\th)\}_{(t,\th)\in\R\times (0,1)}$ be the field defined in (\ref{eq:fieldB}). 
There exists an event $\Omega^*$ of probability $1$ and there is a modification of $\{B(t,\th)\}_{(t,\th)\in\R\times (0,1)}$ denoted by $\{\widetilde{B}(t,\th)\}_{(t,\th)\in\R\times (0,1)}$, such that, for each $\o\in\Omega^*$, the following four results hold.
\begin{itemize}
\item[(i)] The function $(t,\th)\mapsto \widetilde{B}(t,\th,\o)$ is continuous over $\R\times (0,1)$. 
\item[(ii)] For every fixed arbitrarily small real number $\eps >0$ and $(s,\th)\in\R\times(0,1)$, one has
\begin{equation}
\label{eq1:TfieldB}
\limsup_{h\rightarrow 0}\frac{\big |\widetilde{B}(s+h,\th,\o)-\widetilde{B}(s,\th,\o)\big |}{|h|^{\th -\eps}}=0
\end{equation}
and
\begin{equation}
\label{eq2:TfieldB}
\limsup_{h\rightarrow 0}\frac{\big |\widetilde{B}(s+h,\th,\o)-\widetilde{B}(s,\th,\o)\big |}{|h|^{\th +\eps}}=+\infty.
\end{equation}
\item[(iii)] For each fixed $t\in\R$, the function $\th\mapsto \widetilde{B}(t,\th,\o)$ is $C^\infty$ over $(0,1)$; its
derivative, of any order $n\in\Z_+$, at $\th$, is denoted by $\big (\partial_{\th}^{n} \widetilde{B}\big )(t,\th,\o)$.
\item[(iv)] For every fixed $n\in\Z_+$, arbitrarily small real number $\eps>0$ and real numbers $M,a,b$ satisfying $M>0$ and $0<a<b<1$, there exists a constant $C(\o)>0$, only depending on $\o,n,\eps,M,a,b$, such that the inequality, 
\begin{eqnarray}
\label{eq3:TfieldB}
\Big |\big (\partial_{\th}^n \widetilde{B}\big )(t_1,\th_1,\o)-\big (\partial_{\th}^n \widetilde{B}\big )(t_2,\th_2,\o)\Big |\le C(\o)\Big ( |t_1-t_2|^{\max\{\th_1,\th_2\}-\eps}+|\th_1-\th_2|\Big ),\nonumber\\
\end{eqnarray}
holds, for all $(t_1, \th_1)\in [-M,M]\times [a,b]$ and $(t_2, \th_2)\in [-M,M]\times [a,b]$.
\end{itemize}
\end{theorem}

\begin{remark}
\label{rem:AT}
\begin{itemize}
\item The field $\widetilde{B}$ was introduced in \cite{AT03} and Theorem~\ref{theo:fieldB} Parts $(i)$ and $(iii)$ were derived
in the latter article (see \cite{AT03}, pages 463 to 470); notice that in \cite{AT03}, $\widetilde{B}$ was denoted by $B$.
\item A less precise inequality than (\ref{eq3:TfieldB}), was obtained in \cite{AT03}, in the particular case where $n=0$ and $[-M,M]$ is replaced by 
$[0,1]$ (see in \cite{AT03}, Theorem~2.1 and Proposition~2.2 Part $(b)$).
\end{itemize}
\end{remark}

It is worth noticing that a straightforward consequence of Part $(ii)$ of
Theorem~\ref{theo:fieldB}, is the following:
\begin{proposition}
\label{prop:phe-unif-th} For all fixed $\th\in (0,1)$, we denote by
$B_\th =\{B_\th (t)\}_{t\in\R}$ the process $\{\widetilde{B}(t,\th)\}_{t\in\R}$; observe that $B_\th$ is a fBm of Hurst parameter $\th$.
There exists an event $\Omega^*$ of probability $1$, non depending on $s$ and $\th$, such 
that one has, for each $\o\in\Omega^*$ and for all $(s,\th)\in\R\times (0,1)$,
\begin{equation*}
\label{eq:phe-unif-th} 
\a_{B_\th} (s,\o)=\th,
\end{equation*}
where $\a_{B_\th} (s,\o)$ is the \phe at $s$ of the function $t\mapsto B_\th(t,\omega)$. 
\end{proposition}
Observe that the fact that the \phe of the fBm $B_\th$, is equal, almost surely for all $s\in\R$, to the Hurst parameter $\th$, is a classical result (see for example \cite{Xiao97a,Adler81,AJT07}); the novelty 
in Proposition~\ref{prop:phe-unif-th}, is that this equality holds on an event $\Omega^*$ of probability $1$, which does not depend on the Hurst parameter $\th$ (notice that the event $\Omega^*$ also does not depend on $s$).

The proof of Theorem~\ref{theo:fieldB} mainly relies on wavelet techniques, rather similar to those used in \cite{AT03,AJT07}; it is not really the core of the article, this is why it is given in Subsection~\ref{subsec:wav}.\\

From now on, it is important that the reader keeps in his mind the following remark.
\begin{remark}
\label{rem:identBX}
In the remaining of this section as well as in the next section, 
\begin{itemize}
\item the Gaussian field $\{B(t,\th)\}_{(t,\th)\in\R\times (0,1)}$ defined in (\ref{eq:fieldB}), will be always identified with its modification 
$\{\widetilde{B}(t,\th)\}_{(t,\th)\in\R\times (0,1)}$ introduced in Theorem~\ref{theo:fieldB};
\item the mBm $\{X(t)\}_{t\in\R}$ of functional parameter $H$, defined in (\ref{eq:defmBm}), will be always identified with its modification $\{\widetilde{X}(t)\}_{t\in\R}$, defined for every real number $t$ and all $\o\in\Omega$ (the underlying probability space), as $\widetilde{X}(t,\o)=\widetilde{B}(t,H(t),\o)$.
\end{itemize}
\end{remark}

Let us now state the main result of our article.
\begin{theorem}
\label{theo:main} Let $H:\R\rightarrow (0,1)$ be a continuous function, which is nowhere differentiable on some open non-empty interval ${\cal J}\subseteq (0,+\infty)$ and satisfies 
on it, the condition:
$$
\a_H (s)<H(s)<2\a_H (s), \mbox{ for each $s\in {\cal J}$,}
\eqno({\cal A})
$$
where $\a_H (s)$ is the \phe of $H$ at $s$. We denote by $\{\a_X (s)\}_{s\in\R}$, the \phe of $\{X(t)\}_{t\in\R}=\{B(t,H(t))\}_{t\in\R}$, the mBm of functional parameter $H$, and we assume that $\Omega^*$ is the event of 
probability $1$, introduced in Theorem~\ref{theo:fieldB}. Then, the following four results hold.
\begin{itemize}
\item[(i)] For all $\o\in\Omega^*$ and $s\in {\cal J}$, satisfying $(\partial_\th B)(s,H(s),\o)\ne 0$, one has 
$
\a_X (s,\o) =\a_H (s).
$
\item[(ii)] For all $\o\in\Omega^*$ and $s\in {\cal J}$, satisfying $(\partial_\th B)(s,H(s),\o)= 0$, one has 
$
\a_X (s,\o) =H (s).
$
\item[(iii)] There exists $\Omega^{**}\subseteq\Omega^{*}$ an event of probability $1$, such that for all $\o\in\Omega^{**}$,
$$
\dim \big \{s\in {\cal J}\,:\,\,(\partial_\th B)(s,H(s),\o)\ne 0\big\}=1,
$$
where $\dim (\cdot)$ denotes the Hausdorff dimension; in other words,
$$
\dim \big \{s\in {\cal J}\,:\,\,\a_X (s,\o) =\a_H (s)\big\}=1.
$$
\item[(iv)] For each arbitrarily small $\eta>0$, there exists $D\subseteq \Omega^{**}$, an event of (strictly) positive probability, which a priori depends on $\eta$, such that for all $\o\in D$, 
$$
\dim \big \{s\in {\cal J}\,:\,\,(\partial_\th B)(s,H(s),\o)= 0\big\}\ge 1-\eta-\inf_{s\in {\cal J}} H(s)>0;
$$
in other words,
$$
\dim \big \{s\in {\cal J}\,:\,\,\a_X (s,\o) =H (s)\big\}\ge 1-\eta-\inf_{s\in {\cal J}} H(s)>0.
$$
\end{itemize}
\end{theorem}

\begin{remark}
\label{rem:exampleH}
Notice that, when $H$ is a continuous function such that for all $s\in {\cal J}$, one has $H(s)\in [1/3,2/5]$ and $\a_H (s)\in [1/4, 7/24]$, then
condition $({\cal A})$ above, holds. The class of such functions is rather large, namely using methods introduced in \cite{And,DLM,J1,AJT07}, one can 
explicitly construct many of them: for example let $H:\R\rightarrow [1/3,2/5]$ be the $1$-periodic function, defined for 
all $s\in [0,1]$, 
$$
H(s)=\frac{1}{3}+\frac{\left(1-2^{-1/4}\right)}{15}\sum_{j=0}^{+\infty}\sum_{k=0}^{2^j-1} 2^{-j\zeta(k 2^{-j})} T(2^j s-k),
$$
where $\zeta:[0,1]\rightarrow [1/4, 7/24]$ is an arbitrary $C^1$ function satisfying $\zeta(0)=\zeta(1)$ and, where $T:\R\rightarrow [0,1]$ is the function defined,
for each $x\in\R$, as, $T(x)=1-|2x-1|$ if $x\in [0,1]$ and $T(x)=0$ else; by slightly adapting the proof of Proposition~5 in \cite{DLM}, one can show that $\a_H (s)=\zeta(s)\in [1/4, 7/24]$ for any real number $s$.
\end{remark}


The proofs of Parts $(iii)$ and $(iv)$ of Theorem~\ref{theo:main} are postponed to the next section, since they require a specific treatment.
Parts $(i)$ and $(ii)$ will result from Theorem~\ref{theo:fieldB}, let us present the main ideas of their proof, before giving the technical 
details of it. To this end, it is convenient to introduce the following concise notation: let $f$ be a real-valued function defined on neighborhood of $0$, and 
let $\tau\in [0,+\infty)$ be fixed, we assume that the notation:
$$
|f(h)| \asymp |h|^{\tau},
$$
means that for all arbitrarily small $\eps>0$, one has:
$$
\limsup_{h\rightarrow 0} \frac{|f(h)|}{|h|^{\tau-\eps}}=0 \mbox{ and } \limsup_{h\rightarrow 0} \frac{|f(h)|}{|h|^{\tau+\eps}}=+\infty.
$$

\noindent {\sc Heuristic proof of Theorem~\ref{theo:main} Parts $(i)$ and $(ii)$:} For all fixed $\o\in\Omega^*$ and $s\in{\cal J}$, the increment $X(s+h,\o)-X(s,\o)$, of the mBm $\{X(t)\}_{t\in\R}=\{B(t,H(t))\}_{t\in\R}$, can be expressed as:
\begin{equation}
\label{eq1:heurmain}
X(s+h,\o)-X(s,\o)=\big (\Delta B_{H(s)}\big)(s,h,\o)+ R(s,h,\o),
\end{equation}
where
\begin{equation}
\label{eq2:heurmain}
\big (\Delta B_{H(s)}\big)(s,h,\o)=B(s+h,H(s),\o)-B(s,H(s),\o)
\end{equation}
and 
\begin{equation}
\label{eq3:heurmain}
R(s,h,\o)=B(s+h,H(s+h),\o)-B(s+h,H(s),\o).
\end{equation}
Moreover, taking $\th=H(s)$, in (\ref{eq1:TfieldB}) and (\ref{eq2:TfieldB}), one gets that,
\begin{equation}
\label{eq4:heurmain}
\big|\big (\Delta B_{H(s)}\big)(s,h,\o)\big|\asymp |h|^{H(s)}.
\end{equation}
Let us estimate $R(s,h,\o)$. Applying the Mean Value Theorem, it follows that,
\begin{equation}
\label{eq5:heurmain}
R(s,h,\o)=\big(H(s+h)-H(s)\big)\times(\partial_\th B)\big(s+h,\widetilde{\th}(s,h,\o),\o\big),
\end{equation}
where 
\begin{equation}
\label{eq6:heurmain}
\widetilde{\th}(s,h,\o)\in \big (\min\{H(s+h),H(s)\},\max\{H(s+h),H(s)\}\big).
\end{equation}
Next, observe that the definition of the \phe $\a_H (s)$, implies that,
\begin{equation}
\label{eq7:heurmain}
\big|H(s+h)-H(s)\big| \asymp |h|^{\a_H(s)}.
\end{equation} 
Also observe that, in view of (\ref{eq6:heurmain}) and the fact that $(t,\th)\mapsto (\partial_\th B)(t,\th,\o)$ is 
a continuous function, one has that,
\begin{equation}
\label{eq8:heurmain}
(\partial_\th B)\big(s+h,\widetilde{\th}(s,h,\o),\o\big)\xrightarrow[h\rightarrow 0]{} (\partial_\th B)\big(s,H(s),\o\big).
\end{equation}
Next we study two cases: $(\partial_\th B)\big(s,H(s),\o\big)\ne 0$ and $(\partial_\th B)\big(s,H(s),\o\big)=0$.
In the case where $(\partial_\th B)\big(s,H(s),\o\big)\ne 0$; (\ref{eq5:heurmain}),  (\ref{eq7:heurmain}) and (\ref{eq8:heurmain}), imply that 
\begin{equation}
\label{eq9:heurmain}
\big| R(s,h,\o)\big| \asymp |h|^{\a_H(s)}.
\end{equation}
Then putting together, (\ref{eq1:heurmain}), (\ref{eq4:heurmain}), (\ref{eq9:heurmain}) and the inequality $\a_H(s)<H(s)$, one obtains that 
$$
\big|X(s+h,\o)-X(s,\o)\big|\asymp |h|^{\a_H(s)},
$$
which proves that Part $(i)$ of the theorem holds. In the case where $(\partial_\th B)\big(s,H(s),\o\big)= 0$; (\ref{eq3:TfieldB}), (\ref{eq6:heurmain}), 
(\ref{eq7:heurmain}) and the inequality $\a_H(s)<H(s)$, entail that, for all arbitrarily small $\eps>0$,
\begin{equation}
\label{eq10:heurmain}
\big|(\partial_\th B)\big(s+h,\widetilde{\th}(s,h,\o),\o\big)\big| =\mathcal{O}\big(|h|^{\a_H(s)-\eps}\big).
\end{equation}
Then (\ref{eq5:heurmain}), (\ref{eq7:heurmain}) and (\ref{eq10:heurmain}) imply that for all arbitrarily small $\eps>0$,
\begin{equation}
\label{eq11:heurmain}
\big| R(s,h,\o)\big|=\mathcal{O}\big(|h|^{2\a_H(s)-\eps}\big).
\end{equation}
Finally, it follows from (\ref{eq1:heurmain}), (\ref{eq4:heurmain}), (\ref{eq11:heurmain}) and the inequality $H(s)<2\a_H(s)$, that 
$$
\big|X(s+h,\o)-X(s,\o)\big|\asymp |h|^{H(s)},
$$
which proves that Part $(ii)$ of the theorem holds. 
$\Box$

In order to give a rigorous proof of Parts $(i)$ and $(ii)$ of Theorem~\ref{theo:main}, we need some preliminary results. In the remaining of this section, for the sake of simplicity, we assume that ${\cal J}=(0,1)$.

Let us first give two lemmas which, generally
speaking (even in the case where condition~$({\cal A})$ fails to be satisfied), provide upper and lower bounds, for absolute increments
of a typical trajectory of the mBm $X$, in a neighborhood of an
arbitrary fixed point $s\in [0,1]$. 

\begin{lemma}
\label{prop:ub-incremX} For all fixed arbitrarily small $\eps>0$,
$\o\in\Omega^*$ (the event of probability
$1$ introduced in Theorem~\ref{theo:fieldB}) and $s\in [0,1]$, there
is a constant $C(\o)>0$, such that the following inequality, 
\begin{eqnarray*}
\big |X(s+h,\o)-X(s,\o)\big|&\le & C(\o)\left(|h|^{\min\{H(s),2\a_H(s)\}-\eps}+\Big|\big(\partial_\th B\big)(s,H(s),\o)\Big |\times |h|^{\a_H(s)-\eps}\right),
\end{eqnarray*}
holds, for every real number $h$ satisfying $s+h\in [0,1]$.
\end{lemma}

\begin{lemma}
\label{lem1:ex-mBm-sphe} For all fixed arbitrarily small $\eps>0$,
$\o\in\Omega^*$ and $s\in [0,1]$, there exist two constants
$C(\o)>0$ and $C'(\o)>0$ , such that, the inequalities,
\begin{eqnarray}
\label{eq3-bis:ex-mBm-sphe}
&&\big |X(s+h,\o)-X(s,\o)\big|\\
\nonumber
&& \ge \Big|\big(\partial_\th B\big)(s,H(s),\o)\Big |\times \big |H(s+h)-H(s)\big |-C(\o)|h|^{\min\{H(s),2\a_H (s)\}-\eps}\
\end{eqnarray}
and
\begin{eqnarray}
\label{eq3:ex-mBm-sphe}
\big |X(s+h,\o)-X(s,\o)\big| &\ge & \big |B(s+h,H(s),\o)-B(s,H(s),\o)\big |\\
\nonumber
&&-C'(\o)\left(\Big|\big(\partial_\th B\big)(s,H(s),\o)\Big |\times |h|^{\a_H (s)-\eps}+|h|^{\min\{H(s)+\a_H(s),2\a_H (s)\}-\eps}\right).
\end{eqnarray}
hold, for every real number $h$, satisfying $s+h\in [0,1]$.
\end{lemma}

\begin{remark}
\label{rem:lb-pheX} Let $\o\in\Omega^*$ and $s\in (0,1)$ be fixed.
Recall that the \phe at $s$, of the function $t\mapsto X(t,\o)$, is
denoted by $\a_X (s,\o)$. It follows from the previous lemma that:
\begin{itemize}
\item[(i)] When $\big(\partial_\th B\big)(s,H(s),\o)\ne 0$, one has
$$
\a_X (s,\o)\ge \min\{H(s),\a_H(s)\}.
$$
\item[(ii)] When $\big(\partial_\th B\big)(s,H(s),\o)= 0$, one has
$$
\a_X (s,\o)\ge \min\{H(s),2\a_H(s)\}.
$$
\end{itemize}
\end{remark}

Let us now give the proofs of these two lemmas.

\noindent {\sc Proof of Lemma~\ref{prop:ub-incremX}:} First observe that using the first equality in (\ref{eq:defmBm}) as well as the triangle inequality, one has that
\begin{eqnarray}
\label{eq1:ub-incremX}
&& \Big |X(s+h,\o)-X(s,\o)\Big |\\
&& \le \Big | B(s+h,H(s+h),\o)-B(s+h,H(s),\o)\Big |+\Big |B(s+h,H(s),\o)-B(s,H(s),\o)\Big|\nonumber.
\end{eqnarray}
The function $\th\mapsto B(s+h,\th,\o)$ being continuously
differentiable over $(0,1)$ (see Part $(iii)$ of
Theorem~\ref{theo:fieldB}), it follows from the Mean Value Theorem
that there is 
$$
\widetilde{\th}(s,h,\o)\in \big
(\min\{H(s+h),H(s)\},\max\{H(s+h),H(s)\}\big),
$$ 
such that
\begin{equation}
\label{eq2:ub-incremX}
\Big |B(s+h,H(s+h),\o)-B(s+h,H(s),\o)\Big |=\Big|\big(\partial_\th B\big)\big(s+h,\widetilde{\th}(s,h,\o),\o\big)\Big |\times \big |H(s+h)-H(s)\big|.
\end{equation}
Moreover, the triangle inequality implies that
\begin{equation}
\label{eq2-bis:ub-incremX}
\Big|\big(\partial_\th B\big)\big(s+h,\widetilde{\th}(s,h,\o),\o\big)\Big |\le \Big |\big(\partial_\th B\big)\big(s,H(s),\o\big)\Big|
+\Big|\big(\partial_\th B\big)\big(s+h,\widetilde{\th}(s,h,\o),\o\big)-\big(\partial_\th B\big)\big(s,H(s),\o\big)\Big|.
\end{equation}
Part $(iv)$ of Theorem~\ref{theo:fieldB} (in which $\eps$ is replaced by $\eps/2$ and one takes $n=1$, $M=1$, $a=\min_{x\in [0,1]}H(x)$ and $b=\max_{x\in [0,1]}H(x)$), entails that
\begin{eqnarray}
\label{eq3:ub-incremX}
\nonumber
&& \Big|\big(\partial_\th B\big)\big(s+h,\widetilde{\th}(s,h,\o),\o\big)-\big(\partial_\th B\big)\big(s,H(s),\o\big)\Big|\\
&& \le  C_1(\o)\Big (|h|^{\max\{H(s),\widetilde{\th}(s,h,\o)\}-\eps/2}+\big|\widetilde{\th}(s,h,\o)-H(s)\big |\Big )\nonumber\\
&& \le C_1(\o)\Big (|h|^{H(s)-\eps/2}+\big |H(s+h)-H(s)\big|\Big),
\end{eqnarray}
where $C_1(\o)$ is a constant non depending on $s$ and $h$. Putting together (\ref{eq2:ub-incremX}), (\ref{eq2-bis:ub-incremX}),
(\ref{eq3:ub-incremX}) and the inequality $|H(s+h)-H(s)|\le c
|h|^{\a_{H}(s)-\eps/2}$ ($c$ being a constant), one gets that
\begin{eqnarray}
\label{eq3-bis:ub-incremX}
&&\Big |B(s+h,H(s+h),\o)-B(s+h,H(s),\o)\Big |\\
\nonumber
&& \le C_2 (\o) \left (\Big|\big(\partial_\th B\big)(s,H(s),\o)\Big |\times |h|^{\a_H (s)-\eps/2}+|h|^{\min\{H(s)+\a_H (s),2\a_H (s)\}-\eps}\right),
\end{eqnarray}
where $C_2(\o)>0$ is a constant only depending on $\o$, $s$ and $\eps$.
On the other hand, (\ref{eq1:TfieldB}) implies that
\begin{equation}
\label{eq4:ub-incremX}
\Big |B(s+h,H(s),\o)-B(s,H(s),\o)\Big|\le C_3(\o) |h|^{H(s)-\eps},
\end{equation}
where $C_3(\o)>0$ is a constant only depending on $\o$, $s$ and $\eps$. Finally, putting together (\ref{eq1:ub-incremX}), (\ref{eq3-bis:ub-incremX}) and (\ref{eq4:ub-incremX}) one obtains the lemma.
$\Box$
\\

\noindent {\sc Proof of Lemma \ref{lem1:ex-mBm-sphe}:} Let us first show that the inequality (\ref{eq3-bis:ex-mBm-sphe}) is true. It follows from the first equality in (\ref{eq:defmBm}) as well as the triangle inequality, that
\begin{equation}
\label{eq4:ex-mBm-sphe}
\Big |X(s+h,\o)-X(s,\o)\Big |\ge \Big | B(s+h,H(s+h),\o)-B(s+h,H(s),\o)\Big |-\Big |B(s+h,H(s),\o)-B(s,H(s),\o)\Big|.
\end{equation}
Recall that $\Big | B(s+h,H(s+h),\o)-B(s+h,H(s),\o)\Big |$ satisfies the equality (\ref{eq2:ub-incremX}) and that $\Big |B(s+h,H(s),\o)-B(s,H(s),\o)\Big|$ 
satisfies the inequality (\ref{eq4:ub-incremX}); also notice that the triangle
inequality implies that
\begin{equation}
\label{eq5:ex-mBm-sphe}
\Big|\big(\partial_\th B\big)\big(s+h,\widetilde{\th}(s,h,\o),\o\big)\Big |\ge \Big |\big(\partial_\th B\big)\big(s,H(s),\o\big)\Big|
-\Big|\big(\partial_\th B\big)\big(s+h,\widetilde{\th}(s,h,\o),\o\big)-\big(\partial_\th B\big)\big(s,H(s),\o\big)\Big|,
\end{equation}
where $\widetilde{\th}(s,h,\o)$ is as in (\ref{eq2:ub-incremX}).
Putting together (\ref{eq2:ub-incremX}), (\ref{eq5:ex-mBm-sphe}),
(\ref{eq3:ub-incremX}) and the fact that $H$ is a bounded function,
one gets that,
\begin{eqnarray}
\label{eq6:ex-mBm-sphe}
\nonumber && \Big |B(s+h,H(s+h),\o)-B(s+h,H(s),\o)\Big |\\
\nonumber && \ge \left(\Big |\big(\partial_\th B\big)\big(s,H(s),\o\big)\Big|
-\Big|\big(\partial_\th B\big)\big(s+h,\widetilde{\th}(s,h,\o),\o\big)-\big(\partial_\th B\big)\big(s,H(s),\o\big)\Big|\right)\times
\big |H(s+h)-H(s)\big |\\
\nonumber &&\ge \Big|\big(\partial_\th B\big)(s,H(s),\o)\Big |\times\big |H(s+h)-H(s)\big |\\
\nonumber && \hspace{3cm}-C_1 (\o)\Big (|h|^{H(s)-\eps/2}+\big |H(s+h)-H(s)\big|\Big)\times\big |H(s+h)-H(s)\big |\\
&& \ge \Big|\big(\partial_\th B\big)(s,H(s),\o)\Big |\times\big |H(s+h)-H(s)\big|-C_2 (\o)\Big (|h|^{H(s)-\eps/2}+\big
|H(s+h)-H(s)\big|^2\Big),
\end{eqnarray}
where $C_1(\o)>0$ and $C_2(\o)>0$ are two constants only depending on $\o$, $s$ and
$\eps$. Putting together (\ref{eq4:ub-incremX}),
(\ref{eq6:ex-mBm-sphe}), (\ref{eq4:ex-mBm-sphe}) and the inequality $|H(s+h)-H(s)|\le c
|h|^{\a_{H}(s)-\eps/2}$ ($c$ being a constant), one obtains (\ref{eq3-bis:ex-mBm-sphe}). Let us now show that (\ref{eq3:ex-mBm-sphe}) holds. It follows 
from the first equality in (\ref{eq:defmBm}) as well as the triangle inequality that
$$
\Big |X(s+h,\o)-X(s,\o)\Big |\ge \Big |B(s+h,H(s),\o)-B(s,H(s),\o)\Big|-\Big | B(s+h,H(s+h),\o)-B(s+h,H(s),\o)\Big |.
$$
Then combining the latter inequality with (\ref{eq3-bis:ub-incremX}), one gets (\ref{eq3:ex-mBm-sphe}).
$\Box$ \\

\noindent {\sc Rigorous proof of Theorem~\ref{theo:main} Parts $(i)$ and $(ii)$:} The proof can easily be done by making use of condition $({\cal A})$, Remark~\ref{rem:lb-pheX}, Lemma~\ref{lem1:ex-mBm-sphe}, (\ref{eq2:TfieldB}) in which one takes $\th=H(s)$, and the fact 
$$
\limsup_{h\rightarrow 0}\frac{\big|H(s+h)-H(s)\big|}{|h|^{\a_H (s)+\eps}}=+\infty,
$$
for every $\eps>0$.
$\Box$

\section{Hausdorff dimension of the zero-level set of the process $\{(\partial_\th B)(s,H(s))\}_{s\in{\cal J}}$}
\label{sec:levelset} 

The goal of this section is to show that Theorem~\ref{theo:main} Parts $(iii)$ and $(iv)$
 hold. Notice that in all the sequel, we do
not necessarily impose on the continuous functional parameter $H$ of mBm to satisfy condition
(${\cal A}$) (see Theorem~\ref{theo:main}). Also, in all the sequel, we denote by $\{Y(s)\}_{s\in\R}$ the centered Gaussian process $\{(\partial_\th B)(s,H(s))\}_{s\in\R}$, where the centered Gaussian field $\{\big(\partial_\th B\big)(t,\th)\}_{(t,\th)\in\R\times (0,1)}=\{\big(\partial_\th \widetilde{B}\big)(t,\th)\}_{(t,\th)\in\R\times (0,1)}$ has been introduced in Theorem~\ref{theo:fieldB} Part $(iii)$.

Let us first give stochastic integral representations (modifications) of $\{\big(\partial_\th B\big)(t,\th)\}_{(t,\th)\in\R\times (0,1)}$ and $\{Y(s)\}_{s\in\R}$.

\begin{proposition}
\label{prop:rep-partialB}
One has for all $(t,\th)\in \R\times (0,1)$, almost surely,
\begin{equation*}
\label{eq1:rep-partialB}
\big(\partial_\th B\big)(t,\th)=\int_{\R}\Big\{(t-x)_{+}^{\th-1/2}\log\big [(t-x)_+\big ]-(-x)_{+}^{\th-1/2}\log\big [(-x)_+\big ]\Big\}\,dW(x),
\end{equation*}
with the convention that $(y)_{+}^{\th-1/2}\log\big [(y)_+\big ]=0$ for every real numbers $\th\in (0,1)$ and $y\le 0$. As a straightforward consequence, one has for 
all $s\in\R$, almost surely,
\begin{equation}
\label{eq1:rep-Y}
Y(s)=\int_{\R}\Big\{(s-x)_{+}^{H(s)-1/2}\log\big [(s-x)_+\big ]-(-x)_{+}^{H(s)-1/2}\log\big [(-x)_+\big ]\Big\}\,dW(x).
\end{equation}
\end{proposition}

The proof of Proposition~\ref{prop:rep-partialB} is given in Subsection~\ref{subsec:rep-partialB}, since it relies on wavelet techniques similar to those used in Subsection~\ref{subsec:wav}, in order to show that Theorem~\ref{theo:fieldB} holds.

From now on, we assume that $I=[\de_1,\de_2]\subset {\cal J}$ is a compact interval such that $0<\de_2-\de_1<1$ and
\begin{equation}
\label{eq:interv-I}
0< b=\max_{s\in I} H(s)\le \eta+\inf_{s\in {\cal J}} H(s)<1,
\end{equation}
where $\eta>0$ is an arbitrarily small fixed real number; the open interval ${\cal J}$ has been introduced in Theorem~\ref{theo:main} and one has 
that $\de_1>0$ since ${\cal J}\subseteq (0,+\infty)$. Observe that such an interval $I$ exists, since we assume that $H$ is a continuous function over $\R$. The following lemma shows 
that $\Var\big (Y(s)\big)$, $s\in I$, is bounded away from zero. 
\begin{lemma}
\label{lem:lb-varX} There is a constant $c>0$, only depending
on $\de_1$, such that for all $s\in I$, one has
\begin{equation}
\label{eq1:lb-varX}
\Var\big (Y(s)\big)\ge c.
\end{equation}
\end{lemma}

\noindent {\sc Proof of Lemma~\ref{lem:lb-varX}:} It follows from
(\ref{eq1:rep-Y}) and the change of variable $v=s-x$, that
\begin{eqnarray*}
\Var\big (Y(s)\big)&=& \int_{\R}\Big |(s-x)_{+}^{H(s)-1/2}\log\big [(s-x)_+\big]-(-x)_{+}^{H(s)-1/2}\log\big [(-x)_+\big ]\Big |^2\,dx\\
&\ge & \int_{0}^{s}(s-x)^{2H(s)-1}\log^2\big[(s-x)]\,dx\\
&=& \int_{0}^{s}v^{2H(s)-1}\log^2(v)\,dv\\
&\ge & \int_{0}^{\min(\de_1,1)}v\log^2(v)\,dv>0.
\end{eqnarray*}
$\Box$\\

Now we are in position to prove Theorem~\ref{theo:main} Part $(iii)$.

\noindent {\sc Proof of Theorem~\ref{theo:main} Part $(iii)$:} Let us fix $s_0\in I\subset {\cal J}$. Observe that the centered Gaussian random variable 
$Y(s_0)=(\partial_\th B)(s_0,H(s_0))$ has a non-zero standard deviation (see Lemma~\ref{lem:lb-varX}) and thus its probability density function exists; the fact that 
the latter function is strictly positive on the whole real line, implies that $(\partial_\th B)(s_0,H(s_0))\ne 0$ almost surely. Then, noticing that the event $\Omega^*$ (see Theorem~\ref{theo:fieldB}) is of probability $1$, it follows that the 
probability of the event
$$
\Omega^{**}=\big\{\o\in\Omega^{*}:(\partial_\theta B)(s_0,H(s_0),\o)\ne 0\big\},
$$
is equal to $1$. Next, in view of the fact that, for all fixed $\o\in\Omega^{**}\subseteq\Omega^*$, $(\partial_\theta B)(s_0,H(s_0),\o)\ne 0$ and $s\mapsto (\partial_\th B)(s,H(s),\o)$ is a continuous function over the real line (this easily results from Theorem~\ref{theo:fieldB} Part~$(iv)$ as well as from the continuity of $H$), one has that 
$$
\big \{s\in {\cal J}\,:\,\,(\partial_\th B)(s,H(s),\o)\ne 0\big\}={\cal J}\cap\big((\partial_\th B)(\cdot,H(\cdot),\o)\big)^{-1}\big(\R\setminus\{0\}\big),
$$
is a non-empty open subset of $\R$, which implies that its Hausdorff dimension is equal to $1$.
$\Box$\\

The proof of Theorem~\ref{theo:main} Part $(iv)$, is in the same spirit as
that of Theorem 8.4.2 in \cite{Adler81} and Relation (5.9) in \cite{MWX08}; basically, it relies on the following proposition which shows that the process $\{Y(s)\}_{s\in I}$ satisfies the so-called property of one sided strong local nondeterminism. A detailed presentation of the important concept of local nondeterminism and related topics can be found in \cite{Berman73} and in \cite{Xiao06}, for instance.

\begin{proposition}
\label{prop:slnd-Y} For all integer $n\ge 2$ and for each $s^1,\ldots,s^n\in I$ satisfying
\begin{equation}
\label{eq1:slnd-Y}
s^1<\ldots <s^{n},
\end{equation}
one has,
\begin{equation}
\label{eq2:slnd-Y}
\Var\big (Y(s^n)|Y(s^1),\ldots ,Y(s^{n-1})\big )\ge 2^{-1} \big (s^n-s^{n-1}\big )^{2H(s^n)}\log^2\big[(s^n-s^{n-1})\big],
\end{equation}
where $\Var\big (Y(s^n)|Y(s^1),\ldots ,Y(s^{n-1})\big )$ is the conditional variance of $Y(s^n)$ given $Y(s^1),\ldots ,Y(s^{n-1})$.
\end{proposition}

\noindent {\sc Proof of Proposition~\ref{prop:slnd-Y}:} In
view of the definition of $\Var\big (Y(s^n)|Y(s^1),\ldots ,Y(s^{n-1})\big )$ and the Gaussianity of the process $\{Y(s)\}_{s\in I}$, it is sufficient
to show that, for every integer
$n\ge 2$, for all real numbers $a_1,\ldots, a_{n-1}$ and for each
$s^1,\ldots, s^n\in I$ satisfying (\ref{eq1:slnd-Y}), one has
\begin{equation}
\label{eq3:slnd-Y}
\E\left (\Big|Y(s^n)-\sum_{l=1}^{n-1} a_l Y(s^l)\Big|^2\right)\ge 2^{-1} \big (s^n-s^{n-1}\big )^{2H(s^n)}\log^2\big[(s^n-s^{n-1})\big].
\end{equation}
It follows from (\ref{eq1:rep-Y}) and the isometry property of
Wiener integral, that
\begin{eqnarray}
\label{eq4:slnd-Y}
\nonumber
&&\E\left (\Big|Y(s^n)-\sum_{l=1}^{n-1} a_l Y(s^l)\Big|^2\right)\\
&& =\int_{\R}\Big|\Big\{(s^n-x)_{+}^{H(s^n)-1/2}\log\big [(s^n-x)_+\big ]-(-x)_{+}^{H(s^n)-1/2}\log\big [(-x)_+\big ]\Big\}\\
\nonumber
&&\hspace{1.5cm}-\sum_{l=1}^{n-1} a_l \Big\{(s^l-x)_{+}^{H(s^l)-1/2}\log\big [(s^l-x)_+\big ]-(-x)_{+}^{H(s^l)-1/2}\log\big [(-x)_+\big ]\Big\}
\Big|^2\,dx.
\end{eqnarray}
Next, observe that for every $x\in [s^{n-1},s^n]$, one has $-x<0$ and, as a consequence, for all $l\in\{1,\ldots,n\}$,
\begin{equation}
\label{eq4:slnd-Ybis}
(-x)_{+}^{H(s^l)-1/2}\log\big [(-x)_+\big ]=0.
\end{equation}
Also, observe that for every $x\in [s^{n-1},s^n]$ and $l\in\{1,\ldots,n-1\}$, one has, in view of (\ref{eq1:slnd-Y}), that $s^{l}-x\le 0$, therefore,
\begin{equation}
\label{eq4:slnd-Yter}
(s^l-x)_{+}^{H(s^l)-1/2}\log\big [(s^{l}-x)_+\big ]=0.
\end{equation}
Putting together, (\ref{eq4:slnd-Y}), (\ref{eq4:slnd-Ybis}) and (\ref{eq4:slnd-Yter}), one gets that,
\begin{eqnarray}
\label{eq4:slnd-Yquatro}
\nonumber &&\E\left (\Big|Y(s^n)-\sum_{l=1}^{n-1} a_l Y(s^l)\Big|^2\right)\\
\nonumber && \ge \int_{s^{n-1}}^{s^n}(s^n-x)^{2H(s^n)-1}\log^2 (s^n-x)\,dx\\
&& =\int_{0}^{s^n-s^{n-1}}(s^n-s^{n-1}-x)^{2H(s^n)-1}\log^2 (s^n-s^{n-1}-x)\,dx.
\end{eqnarray}
Next, setting in the last integral, $v=x/(s^n-s^{n-1})$ and using the fact that $0<s^n-s^{n-1}\le \de_2-\de_1<1$, one obtains that
\begin{eqnarray}
\label{eq5:slnd-Y}
\nonumber
&& \int_{0}^{s^n-s^{n-1}}(s^n-s^{n-1}-x)^{2H(s^n)-1}\log^2 (s^n-s^{n-1}-x)\,dx\\
\nonumber
&& =(s^n-s^{n-1})\int_0 ^1 \Big [(s^n-s^{n-1})-(s^n-s^{n-1})v\Big]^{2H(s^n)-1}\log^2 \Big [(s^n-s^{n-1})-(s^n-s^{n-1})v\Big]\,dv\\
\nonumber
&&=(s^n-s^{n-1})^{2H(s^n)}\int_0 ^1 (1-v)^{2H(s^n)-1}\Big (\log\big[(1-v)^{-1}\big]+\log\big[ (s^n-s^{n-1})^{-1}\big]\Big)^2\,dv\\
\nonumber
&& \ge (s^n-s^{n-1})^{2H(s^n)}\log^2\big[(s^n-s^{n-1})\big]\int_0 ^1
(1-v)\,dv\\
&& = 2^{-1} (s^n-s^{n-1})^{2H(s^n)}\log^2\big[(s^n-s^{n-1})\big].
\end{eqnarray}
Finally combining (\ref{eq4:slnd-Yquatro}) with (\ref{eq5:slnd-Y}), it follows that 
(\ref{eq3:slnd-Y}) holds. 
$\Box$\\

Now it is convenient to make a few recalls concerning the so-called
(deterministic) measures of finite $\ga$ energy, more information
about them can be found in \cite{Landkof72}. In the sequel, we
always assume that $\ga\in (0,1)$. A measure $\mu$ defined on the
Borel sets of $\R$ is said to be of finite $\ga$ energy, if the
integral
\begin{equation}
\label{eq1:meas-fi}
{\cal I}_\ga (\mu)=\int_{\R}\int_{\R} |s-t|^{-\ga}\,d\mu(s)d\mu(t),
\end{equation}
which is usually called the $\ga$ energy of $\mu$, exists and is finite. ${\cal M}_\ga$, the class of these measures, forms a Hilbert space equipped with the inner product
$$
(\mu,\nu)_\ga= \int_{\R}\int_{\R} |s-t|^{-\ga}\,d\mu(s)d\nu(t);
$$
the corresponding norm is denoted by $\|\cdot\|_\ga$. Moreover, ${\cal M}^+_\ga$, the subset of the positive measures of ${\cal M}_\ga$, is a complete metric space, for the metric
\begin{equation}
\label{eq2:meas-fi}
\|\mu-\nu\|_\ga=\sqrt{\int_{\R}\int_{\R} |s-t|^{-\ga}\,d(\mu-\nu)(s)d(\mu-\nu)(t)}=\sqrt{{\cal I}_\ga (\mu-\nu)}.
\end{equation}
One of the main interests of the positive measures of finite $\ga$ energy
comes from the following lemma which is a straightforward
consequence of the Frostman Theorem, the latter theorem is
presented in e.g. \cite{kahane85} pages 132 and 133 (see also
\cite{Fa1}).
\begin{lemma}
\label{lem:Frostman}
Let $K$ be a compact subset of $\R$. If $K$ carries a positive non-vanishing measure of finite $\ga$ energy (i.e. if there is a positive non-vanishing measure of finite $\ga$ energy whose support is contained in $K$), then the Hausdorff dimension of $K$ is greater than or equal to $\ga$.
\end{lemma}

We are now in position to prove Theorem~\ref{theo:main} Part $(iv)$.

\noindent {\sc Proof of Theorem~\ref{theo:main} Part $(iv)$:} For all
$\o\in\Omega^{*}$ (the event of probability $1$ introduced in Theorem~\ref{theo:fieldB}), we set 
$$
{\cal L}_Y(\o)=\big \{s\in I\,:\,\,(\partial_\th B)(s,H(s),\o)= 0\big\}.
$$
Recall that $I=[\de_1,\de_2]\subset {\cal J}$ is a compact interval such that $0<\de_2-\de_1<1$ and (\ref{eq:interv-I}) holds. Also, recall 
that the process $\{(\partial_\th B)(s,H(s))\}_{s\in I}$ is denoted by $\{Y(s)\}_{s\in I}$.

We will show that for all $\ga<1-b$ there is ${\cal D}_{\ga}\subseteq
\Omega^{*}$ an event, which a priori depends on $\ga$ but whose
probability is bigger than a strictly positive constant non
depending on $\ga$, such that for all $\o\in {\cal D}_{\ga}$, the set ${\cal
L}_{Y}(\o)$ carries a positive non-vanishing deterministic
measure $\mu(\cdot,\o)$, whose $\ga$ energy is finite. Roughly speaking, the idea for obtaining
$\mu(\cdot,\o)$ is somehow similar to the one which consists in constructing a Dirac
measure as a limit of Gaussian measures; namely, 
$\mu(\cdot,\o)$ will be the limit, in the sense of the
norm $\|\cdot\|_\ga$, of some of the positive measures $\mu_n(\cdot,\o)$, $n\in\N$, 
defined for each Borel subset $A$ of $\R$, as,
\begin{equation}
\label{eq4:lsY}
\mu_n(A,\o)=\int_{A\cap I} \Phi_n (t,\o)\,dt,
\end{equation}
where for all $t\in I$,
\begin{equation}
\label{eq5:lsY}
\Phi_n (t,\o)=\big (2\pi n\big)^{1/2}\exp\left [-\frac{nY(t,\o)^2}{2}\right].
\end{equation}
Notice that (\ref{eq5:lsY}) and the Fourier inversion formula, imply that for all $t\in I$ and $\o\in\Omega^{*}$,
\begin{equation}
\label{eq5bis:lsY}
\Phi_n (t,\o)=\int_{\R}\exp\left [-\frac{\xi^2}{2n}+i\xi Y(t,\o)\right]\,d\xi.
\end{equation}
It is clear that
\begin{equation}
\label{eq5ter:lsY}
\Sp \,\mu_n(\cdot,\o)=I.
\end{equation}
Moreover, $\mu_n(\cdot,\o)$ is of finite $\ga$ energy for any $\ga\in (0,1)$. Indeed, in view of (\ref{eq1:meas-fi}) and (\ref{eq4:lsY}), ${\cal I}_\ga\big(\mu_n(\cdot,\o)\big)$ can be expressed as,
\begin{equation}
\label{eq6:lsY}
{\cal I}_\ga\big(\mu_n(\cdot,\o)\big)=\int_I\int_I |s-t|^{-\ga}\Phi_n (s,\o)\Phi_n (t,\o)\,ds\,dt;
\end{equation}
then, (\ref{eq5:lsY}) implies that
$$
{\cal I}_\ga\big(\mu_n(\cdot,\o)\big)\le 2\pi n\int_I\int_I |s-t|^{-\ga}\,ds\,dt<\infty.
$$
Observe that one can more generally, show in the same way, that for
all integers $n\ge 1$ and $m\ge 1$,
$$
\int_I\int_I |s-t|^{-\ga}\Phi_n (s,\o)\Phi_m (t,\o)\,ds\,dt<\infty.
$$
Let us now construct a subsequence $l\mapsto n_l$ satisfying the
following property: for all $\ga<1-b$, there is $\check{\Omega}_{\ga}\subseteq
\Omega^{*}$ an event of probability $1$, which a priori depends on
$\ga$, such that for each $\o\in \check{\Omega}_{\ga}$, $\big
(\mu_{n_l}(\cdot,\o)\big)_l$ is a Cauchy sequence, in the sense of
the norm $\|\cdot\|_\ga$. To this end, one needs to give, for all
integers $n\ge 1$ and $p\ge 0$, a convenient upper bound of the
quantity $\E\big (\|\mu_{n+p}-\mu_n\|_{\ga}^2\big )$. By using (\ref{eq2:meas-fi}), (\ref{eq6:lsY}), Fubini Theorem and
(\ref{eq5bis:lsY}), one gets that,
\begin{eqnarray}
\label{eq7:lsY}
\nonumber
&&\E\big (\|\mu_{n+p}-\mu_n\|_{\ga}^2\big )\\
&&=\int_I \int_I |s-t|^{-\ga}\E\Big\{\big(\Phi_{n+p}(s)-\Phi_n (s)\big )\big(\Phi_{n+p}(t)-\Phi_n (t)\big )\Big\}\,ds\,dt\\
\nonumber
&&=\int_I \int_I \int_\R \int_\R |s-t|^{-\ga}\left (\exp\Big [-\frac{\xi^2}{2(n+p)}\Big]-\exp\Big [-\frac{\xi^2}{2n}\Big]\right )
\left (\exp\Big [-\frac{\e^2}{2(n+p)}\Big]-\exp\Big [-\frac{\e^2}{2n}\Big]\right )\\
\nonumber
&& \hspace{5cm}\times\E\left [\exp\Big(i\big (\xi Y(s)+\e Y(t)\big)\Big)\right]\,d\xi\,d\e\,ds\,dt.
\end{eqnarray}
Moreover, in view of the fact that $(\xi,\e)\mapsto \E\left [\exp\Big(i\big (\xi Y(s)+\e Y(t)\big)\Big)\right]$ is the characteristic function of the centered Gaussian vector
$
\left (
\begin{array}{c}
Y(s)\\
Y(t)
\end{array}
\right ),
$
 one has that
\begin{equation}
\label{eq8:lsY}
\E\left [\exp\Big(i\big (\xi Y(s)+\e Y(t)\big)\Big)\right]=\exp\left\{-\frac{1}{2}
\left (
\begin{array}{c}
\xi\\
\e
\end{array}
\right )^t
\Ga_Y (s,t)
\left (
\begin{array}{c}
\xi\\
\e
\end{array}
\right )\right\},
\end{equation}
where
$
\left (
\begin{array}{c}
\xi\\
\e
\end{array}
\right )^t
$
is the transpose of
$
\left (
\begin{array}{c}
\xi\\
\e
\end{array}
\right )
$
and where
$
\Ga_Y (s,t)=
\left (
\begin{array}{ll}
\E\big [Y(s)^2\big] & \E\big [Y(s)Y(t)\big ] \\
\E\big [Y(s)Y(t)\big ] & \E\big [Y(t)^2\big ]
\end{array}
\right )
$ is the covariance matrix of $
\left (
\begin{array}{c}
Y(s)\\
Y(t)
\end{array}
\right ).
$
Also, observe that for all integers $n\ge 1$, $p\ge 0$ and real number $\xi$,
\begin{eqnarray}
\label{eq9:lsY}
\nonumber
0\le \exp\Big [-\frac{\xi^2}{2(n+p)}\Big]-\exp\Big [-\frac{\xi^2}{2n}\Big]&\le &
\exp\Big [-\frac{\xi^2}{2(n+p)}\Big]\left(1-\exp\Big [-\frac{p}{2n(n+p)}\xi^2\Big ]\right)\\
&\le & 1-\exp\Big [-\frac{1}{2n}\xi^2\Big ].
\end{eqnarray}
Putting together, (\ref{eq7:lsY}), (\ref{eq8:lsY}) and (\ref{eq9:lsY}), one obtains that for all integers $n\ge 1$, $p\ge 0$
\begin{equation}
\label{eq10:lsY}
0\le \E\big (\|\mu_{n+p}-\mu_n\|_{\ga}^2\big )\le U_n,
\end{equation}
where
\begin{eqnarray}
\label{eq11:lsY}
&& U_n=\int_I \int_I \int_\R \int_\R |s-t|^{-\ga}\left (1-\exp\Big [-\frac{1}{2n}\xi^2\Big ]\right )
\left (1-\exp\Big [-\frac{1}{2n}\e^2\Big ]\right )\\
\nonumber
&& \hspace{5cm}\times\exp\left\{-\frac{1}{2}
\left (
\begin{array}{c}
\xi\\
\e
\end{array}
\right )^t
\Ga_Y (s,t)
\left (
\begin{array}{c}
\xi\\
\e
\end{array}
\right )\right\}\,d\xi\,d\e\,ds\,dt.
\end{eqnarray}
Let us now show that
\begin{equation}
\label{eq12:lsY}
\lim_{n\rightarrow +\infty} U_n=0.
\end{equation}
The equality (\ref{eq12:lsY}) results from the dominated convergence Theorem. Indeed, for almost all $(s,t,\xi,\e)\in I^2\times\R^2$, one has
\begin{eqnarray}
\label{eq12bis:lsY}
&&\lim_{n\rightarrow +\infty}|s-t|^{-\ga}\left (1-\exp\Big [-\frac{1}{2n}\xi^2\Big ]\right )
\left (1-\exp\Big [-\frac{1}{2n}\e^2\Big ]\right )\\
\nonumber
&& \hspace{5cm}\times\exp\left\{-\frac{1}{2}
\left (
\begin{array}{c}
\xi\\
\e
\end{array}
\right )^t
\Ga_Y (s,t)
\left (
\begin{array}{c}
\xi\\
\e
\end{array}
\right )\right\}=0.
\end{eqnarray}
Moreover, using the equalities
\begin{equation}
\label{eq13:lsY}
\int_\R \int_\R \exp\left\{-\frac{1}{2}
\left (
\begin{array}{c}
\xi\\
\e
\end{array}
\right )^t
\Ga_Y (s,t)
\left (
\begin{array}{c}
\xi\\
\e
\end{array}
\right )\right\}\,d\xi\,d\e =2\pi\left(\det \big ( \Ga_Y (s,t)\big )\right )^{-1/2},
\end{equation}
and
\begin{equation}
\label{eq14:lsY}
\det \big ( \Ga_Y (t,s)\big )=\Var \big (Y(s)\big )\times \Var \big( Y(t)|Y(s)\big ),
\end{equation}
and using the fact that $I=[\de_1,\de_2]$, Lemma~\ref{lem:lb-varX}, Proposition~\ref{prop:slnd-Y} and (\ref{eq:interv-I}), one has
\begin{eqnarray}
\label{eq15:lsY}
\nonumber
&&\int_I \int_I \int_\R \int_\R |s-t|^{-\ga}\exp\left\{-\frac{1}{2}
\left (
\begin{array}{c}
\xi\\
\e
\end{array}
\right )^t
\Ga_Y (s,t)
\left (
\begin{array}{c}
\xi\\
\e
\end{array}
\right )\right\}\,d\xi\,d\e\,ds\,dt\\
\nonumber
&&=4\pi\int_{\de_1}^{\de_2}\int_{\de_1}^{t}(t-s)^{-\ga}\left (\det \big ( \Ga_Y (s,t)\big )\right)^{-1/2}\,ds\,dt\\
\nonumber
&&=4\pi\int_{\de_1}^{\de_2}\int_{\de_1}^{t}(t-s)^{-\ga}\left (\Var \big (Y(s)\big )\times \Var \big( Y(t)|Y(s)\big )\right )^{-1/2}\,ds\,dt\\
&& \le c_1 \int_{\de_1}^{\de_2}\int_{\de_1}^{t}(t-s)^{-\ga-b}\big|\log(t-s)\big|^{-1}\,ds\,dt<\infty,
\end{eqnarray}
where $c_1$ is a finite constant only depending on $\de_1$; observe that the last integral is finite since $\ga<1-b$.
Relations (\ref{eq12bis:lsY}) and (\ref{eq15:lsY}), allow us to use
the dominated convergence Theorem and to obtain Relation
(\ref{eq12:lsY}). Next, it follows from (\ref{eq12:lsY}), that there is
an increasing subsequence $l\mapsto n_l$ such that for all $l$ one
has
$$
U_{n_l}\le 2^{-l}.
$$
Then setting in (\ref{eq10:lsY}), $n=n_l$ and $p=n_{l+1}-n_l$, and using Cauchy-Schwarz inequality, one obtains that
$$
\E\big (\|\mu_{n_{l+1}}-\mu_{n_l}\|_{\ga}\big )\le 2^{-l/2}
$$
and consequently that
$$
\E\left(\sum_{l=0}^{+\infty} \|\mu_{n_{l+1}}-\mu_{n_l}\|_{\ga}\right)<\infty.
$$
This implies that there exists an event $\check{\Omega}_{\ga}\subseteq\Omega^{*}$ of probability $1$ such that for all $\o\in \check{\Omega}_{\ga}$,
$$
\sum_{l=0}^{+\infty} \|\mu_{n_{l+1}}(\cdot,\o)-\mu_{n_l}(\cdot,\o)\|_{\ga}<\infty.
$$
Therefore $\big(\mu_{n_l}(\cdot,\o)\big)_l$ is a Cauchy sequence for the norm $\|\cdot\|_\ga$ and, as a consequence, it converges to a positive measure $\mu(\cdot,\o)$ of finite $\ga$ energy. In view of (\ref{eq5ter:lsY}), one clearly has that $\Sp\,\mu(\cdot,\o)\subseteq I$, let us show that one even has,
\begin{equation}
\label{eq16:lsY}
\Sp\,\mu(\cdot,\o)\subseteq {\cal L}_{Y}(\o).
\end{equation}
Let $g$ be a bounded continuous function on the real line which vanishes on a neighborhood of ${\cal L}_{Y}(\o)$ and let $K$ be the compact set defined as $K=I\cap (\Sp\,g)$. Observe that, in view of the definition of ${\cal L}_{Y}(\o)$ and of the continuity of the function $t\mapsto Y(t,\o)^2$, there is a constant $C_2 (\o)>0$ such that for all $t\in K$, $Y(t,\o)^2\ge C_2(\o)$. Therefore, using (\ref{eq4:lsY}) and (\ref{eq5:lsY}), one has for all~$l$,
$$
\Big |\int_{\R}g(t)\,d\mu_{n_l}(t,\o)\Big|=(2\pi n_l )^{1/2}\Big |\int_K g(t)e^{-n_l Y(t,\o)^2/2}\,dt\Big |\le (2\pi n_l )^{1/2}e^{-n_l C_{2}(\o)/2}\max_{t\in K} |g(t)|,
$$
and consequently that
$$
\int_{\R}g(t)\,d\mu(t,\o)=\lim_{l\rightarrow+\infty} \int_{\R}g(t)\,d\mu_{n_l}(t,\o)=0;
$$
thus one obtains (\ref{eq16:lsY}).
Let us now show that there are two constants $c_3>0$ and $c_4>0$, which do not depend on $\ga$, and that there exists an event ${\cal D}_{\ga}\subseteq\check{\Omega}_{\ga}$, a priori depending on $\ga$, which satisfies
\begin{equation}
\label{eq17:lsY}
\P({\cal D}_{\ga})\ge c_3
\end{equation}
and for all $\o\in {\cal D}_{\ga}$,
\begin{equation}
\label{eq18:lsY}
\mu ( {\cal L}_{Y}(\o),\o)=\mu(I,\o)\ge c_4.
\end{equation}
To this end, we will use the following lemma whose proof is elementary (see e.g. \cite{kahane85} page 8).

\begin{lemma}
\label{lem:lbprob}
Let $X$ be a real-valued nonnegative random variable with a finite non-vanishing second moment. Then one has for all $\la\in (0,1)$,
\begin{equation}
\label{eq:lbprob}
P\big\{X\ge \la \E(X)\big\}\ge (1-\la)^2\frac{\E^2(X)}{\E(X^2)}.
\end{equation}
\end{lemma}

It follows from (\ref{eq4:lsY}), (\ref{eq5bis:lsY}), Fubini Theorem and the fact that $\xi\mapsto \E[e^{i\xi Y(t)}]$ is the characteristic function of the centered Gaussian random variable $Y(t)$,
that
$$
\E\big [\mu_{n_l}(I)\big ]=\int_I \int_\R e^{-\xi^2/2n_l}\E \big [ e^{i\xi Y(t)}\big]\,d\xi\,dt=\int_I \int_\R e^{-\xi^2/2n_l} e^{-\si_Y (t)^2\xi^2/2}\,d\xi\,dt,
$$
where $\si_Y (t)$ is the standard deviation of $Y(t)$. Then, using the dominated convergence Theorem and Lemma~\ref{lem:lb-varX}, one has that
\begin{equation}
\label{eq19:lsY}
\lim_{l\rightarrow +\infty} \E\big [\mu_{n_l}(I)\big ]=\sqrt{2\pi}\int_I \si_Y(t)^{-1}\,dt=c_5>0.
\end{equation}
On the other hand, (\ref{eq4:lsY}), (\ref{eq5bis:lsY}), Fubini Theorem, and (\ref{eq8:lsY}) imply that
\begin{eqnarray*}
&&
\E\big [\mu_{n_l}(I)^2\big ]\\
&& =\int_I \int_I\int_\R \int_\R e^{-\frac{\xi^2+\e^2}{2n_l}}\E\big [e^{i(\xi Y(s)+\e Y(t))}\big ]\,d\xi\,d\e\,ds\,dt\\
&& =\int_I \int_I\int_\R \int_\R e^{-\frac{\xi^2+\e^2}{2n_l}}\exp\left\{-\frac{1}{2}
\left (
\begin{array}{c}
\xi\\
\e
\end{array}
\right )^t
\Ga_Y (s,t)
\left (
\begin{array}{c}
\xi\\
\e
\end{array}
\right )\right\}\,d\xi\,d\e\,ds\,dt.
\end{eqnarray*}
Then using the dominated convergence Theorem, (\ref{eq13:lsY}), (\ref{eq14:lsY}), the equality $I=[\de_1,\de_2]$, Lemma~\ref{lem:lb-varX}, Proposition~\ref{prop:slnd-Y}
and (\ref{eq:interv-I}), one obtains that
\begin{eqnarray}
\label{eq20:lsY}
0<c_6=\lim_{l\rightarrow +\infty} \E\big [\mu_{n_l}(I)^2\big ] &=& 2\pi \int_I \int_I \big (\det\, \Ga_Y(s,t)\big )^{-1/2}\,ds\, dt\\
\nonumber
&\le c_1 & \int_{\de_1}^{\de_2} \int_{\de_1}^t (t-s)^{-b}\big|\log(t-s)\big|^{-1}\,ds\,dt<\infty,
\end{eqnarray}
where $c_1$ is the finite constant already introduced in (\ref{eq15:lsY}) and where the last integral is finite since $b<1$.
Let $\la_0\in (0,1)$ be such that
$$
\P\big\{\mu\big (I\big )=\la_0 c_5 \big\}=0.
$$
By using (\ref{eq19:lsY}) as well as the fact that one has, almost
surely (more precisely, on the event~$\check{\Omega}_\ga$),
$$
\mu\big(I\big)=\lim_{l\rightarrow +\infty}\mu_{n_l}(I),
$$
one obtains that
$$
\P\left\{\mu\big (I\big )\ge \la_0 c_5 \right\}=\lim_{l\rightarrow +\infty} \P\left\{\mu_{n_l}\big (I\big )
\ge \la_0\E\big [\mu_{n_l}(I)\big ]\right\}.
$$
Then (\ref{eq19:lsY}), (\ref{eq20:lsY}) and Lemma~\ref{lem:lbprob} 
imply that
$$
\P\left\{\mu\big (I)\ge \la_0 c_5 \right\}\ge (1-\la_0)^2\frac{c_{5}^2}{c_6}.
$$
Therefore, setting $c_3=(1-\la_0)^2\frac{c_{5}^2}{c_6}$, $c_4=\la_0 c_5$ and ${\cal D}_{\ga}=\{\o\in \check{\Omega}_\ga :\mu(I,\o)\ge c_4\}$, one gets (\ref{eq17:lsY}) and (\ref{eq18:lsY}). Next it follows from Lemma~\ref{lem:Frostman} that, for all real number $\ga$, satisfying $\ga<1-b$ and for all $\o\in {\cal D}_{\ga}$, one has
\begin{equation}
\label{eq22:lsY}
\dim {\cal L}_Y(\o)\ge \ga.
\end{equation}
For every integer $m>(1-b)^{-1}$, let $D_{m}$ be the event defined as
$$
D_{m}=\bigcup_{n=m}^{+\infty} {\cal D}_{1-b-n^{-1}}.
$$
It is clear that for all $m$, $D_{m+1}\subseteq D_m$. Moreover, (\ref{eq17:lsY}) implies that
$
\P(D_{m})\ge c_3
$
and (\ref{eq22:lsY}) that for all $\o\in D_{m}$,
$
\dim {\cal L}_Y(\o)\ge 1-b-m^{-1}.
$
Therefore taking
$$
D=\left(\bigcap_{m>(1-b)^{-1}}^{+\infty} D_{m}\right)\cap\Omega^{**},
$$
where $\Omega^{**}$ is the event of probability~$1$ introduced in Theorem~\ref{theo:main} Part~$(iii)$, one obtains Theorem~\ref{theo:main} Part $(iv)$.
$\Box$

\section{Appendix}
\label{sec:appendice}


\subsection{Proof of Theorem~\ref{theo:expdet}}
\label{subsec:expdet}

\noindent {\it Proof of Theorem~\ref{theo:expdet} Part $(i)$:} First observe that by using the same method as in the proofs of Lemma~3.5 in \cite{AT03} 
and Proposition~3.6 in \cite{AJT07}, one can show that there exists a deterministic quantity $\widetilde{b}_X (J)\in [0,1]$, such that,
\begin{equation}
\label{eq3:expdet}
\P\big\{\beta_{X}(J)= \widetilde{b}_X (J)\big\}=1.
\end{equation}
Let us prove that $\widetilde{b}_X (J)=b_X(J)$. It follows from (\ref{eq1:expdet}), (\ref{eq:unifhold}), (\ref{eq3:expdet}) and Lemma~2.3 in \cite{AT03}, that 
$\widetilde{b}_X (J)\ge b_X(J)$. In view of the latter inequality and the fact that $b_X (J)\ge 0$, it is clear that one has 
$\widetilde{b}_X (J)=b_X(J)$ when $\widetilde{b}_X (J)=0$. So from now on, we assume that $\widetilde{b}_X (J)\in (0,1]$. Let 
then $\la$ be an arbitrary deterministic real number belonging to the open interval $(0,\widetilde{b}_X (J))$. Relations (\ref{eq3:expdet}) and (\ref{eq:unifhold}),
imply that 
\begin{equation}
\label{eq4:expdet}
\sup_{t',t''\in J}\frac{|X(t',\o)-X(t'',\o)|}{|t'-t''|^\la}<\infty, \mbox{ almost surely.}
\end{equation}
Next (\ref{eq4:expdet}) and the Gaussianity of the process $\{X(t)\}_{t\in J}$, entail (see \cite{LedTal}) that 
$$
\sup_{t',t''\in J}\frac{\E|X(t',\o)-X(t'',\o)|^2}{|t'-t''|^{2\la}}\le \E\left [\sup_{t',t''\in J}\frac{|X(t',\o)-X(t'',\o)|^2}{|t'-t''|^{2\la}}\right]<\infty
$$
and, as a consequence (see (\ref{eq1:expdet})), that $\la\le b_X(J)$. One gets, from the latter inequality, that $\widetilde{b}_X (J)\le b_X(J)$, since $\la\in(0,\widetilde{b}_X (J))$ is arbitrary.
$\Box$\\

\noindent {\it Proof of Theorem~\ref{theo:expdet} Part $(ii)$:} First observe that, assuming that $J_1$ and $J_2$ are two arbitrary non-degenerate compact intervals satisfying $J_1\subseteq J_2$, one has, in view of (\ref{eq1:expdet}) and (\ref{eq:unifhold}), that,
$$
b_X (J_1)\ge b_X (J_2)\mbox{ and } \beta_X (J_1,\omega)\ge \beta_X (J_2,\omega) \mbox{ for all $\o$.}
$$
Therefore (\ref{eq2:expdet}) and (\ref{eq:lochold1}), imply that,
\begin{equation}
\label{eq:lochold4}
\widetilde{a}_X(s)=\sup\left\{b_X\big([u,v]\big)\,:\,\,u,v\in\Q \mbox{ and } s\in (u,v)\right\}
\end{equation}
and, for each $\o$,
\begin{equation}
\label{eq:lochold}
\widetilde{\a}_X(s,\omega)=\sup\left\{\beta_X\big([u,v],\o\big)\,:\,\,u,v\in\Q \mbox{ and } s\in (u,v)\right\}, 
\end{equation}
where $\Q$ denotes the set of the rational numbers. On the other hand, (\ref{eq:unif-holder}) and the fact that $\Q$ is a countable set, entail that,
\begin{equation}
\label{eq:unif-holder2}
\P\left\{\beta_X\big( [u,v]\big)=b_X\big([u,v]\big)\,:\,\, \mbox{for all $u,v\in\Q$ such that $u<v$}\right\}=1. 
\end{equation}
Putting together, (\ref{eq:lochold4}), (\ref{eq:lochold}) and (\ref{eq:unif-holder2}), one obtains (\ref{eq:stongphl}).
$\Box$

\subsection{Proof of Theorem~\ref{theo:fieldB}}
\label{subsec:wav}
As we have already mentioned (see Remark~\ref{rem:AT}), Parts $(i)$ and $(iii)$ of Theorem~\ref{theo:fieldB} have been obtained in \cite{AT03}, however, we need several ingredients of their proofs, in order to derive the other parts of the theorem. This is the reason why these proofs will be recalled in the sequel.

The modification $\{\widetilde{B}(t,\th)\}_{(t,\th)\in\R\times (0,1)}$  of the field $\{B(t,\th)\}_{(t,\th)\in\R\times (0,1)}$ (see (\ref{eq:fieldB})), 
will be defined as a random wavelet type series. Let us first introduce some notations related to wavelets.
\begin{itemize}
\item We denote by $\psi$ a
Lemari\'e-Meyer real-valued mother wavelet \cite{LM86,meyer92,daubechies92}.
Recall that it satisfies the following three nice properties:
\begin{itemize}
\item[(a)] $\psi$ belongs to the Schwartz class $S(\R)$; which means that $\psi$ is a $C^\infty$ function and decreases at infinity, as well as all its derivatives of any order, faster than any polynomial.
\item[(b)] The support of $\widehat{\psi}$, the Fourier transform of $\psi$, is contained in the ring $\big\{\xi\in\R\,:\,\,\frac{2\pi}{3}\le |\xi|\le \frac{8\pi}{3}\}$; throughout this subsection, the Fourier transform of a function $f\in L^1(\R)$ is definied for every real number $\xi$, as $\widehat{f}(\xi)=\int_\R e^{-i\xi x}f(x)\,dx$.
Sometime, we denote $\widehat{f}$ by ${\cal F}(f)$. 
 
\item[(c)] The collection of the functions:
$$
{\cal ML}=\big\{2^{j/2}\psi (2^j\cdot-k)\,:\,\,(j,k)\in\Z^2\big\},
$$
forms an orthonormal basis of the Hilbert space $L^2(\R)$.
\end{itemize}
Observe that $(a)$ and $(b)$ imply that $\psi$ belongs to Lizorkin
space (see e.g. page 148 in \cite{SKM93} for a definition of this
space).
\item We denote by $\Psi$ the real-valued function defined for every
$(y,\th)\in\R\times (0,1)$ as,
\begin{equation}
\label{eq1:pfpsi}
\Psi(y,\th)=\int_{\R} (y-x)_{+}^{\th-1/2}\psi (x)\,dx.
\end{equation}
Observe that for all fixed $\th\in (0,1)$, the Fourier transform ${\cal F}\big (\Psi(\cdot,\th)\big)$ of the function $t\mapsto \Psi(t,\th)$, satisfies, for all real number $\xi\neq 0$,
\begin{equation}
\label{eq:Fourier-pfpsi}
{\cal F}\big (\Psi(\cdot,\th)\big)(\xi)=\Ga(\th+1/2) e^{-i\,\sgn(\xi) (\th+1/2)\frac{\pi}{2}}
\frac{\widehat{\psi}(\xi)}{|\xi|^{\th+1/2}},
\end{equation}
where $\Ga$ is the usual Gamma function, defined for every real number
$z>0$, as,
$$
 \Ga (z)=\int_{0}^{+\infty} x^{z-1} e^{-x}\,dx.
$$
Relation (\ref{eq:Fourier-pfpsi}) comes from the fact that for each
fixed $\th\in (0,1)$, the function 
$
\frac{\Psi (\cdot,\th)}{\Ga(\th+1/2)},
$
is the left-sided fractional primitive of $\psi$ of order $\th+1/2$,
we refer to Chapter 2 of \cite{SKM93} for its proof. It is worth noticing that by using (\ref{eq:Fourier-pfpsi}) and a method quite similar to that which 
allowed to obtain Lemma 2.1 in \cite{AT03} and Lemma 2.4 in \cite{ASX}, one can show that $\Psi$ is $C^\infty$ over $\R\times (0,1)$ and also, that it is, as well as  its partial derivatives of any order, well-localized in the variable $y\in\R$, uniformly
in the variable $\th\in (0,1)$; in other words, for every nonnegative integers $m$ and $n$, one has
\begin{equation}
\label{eq1:prop-fdpsi} \sup\left\{\big(2+|y|\big)^2\big|(\partial_{y}^m \partial_{\th}^n \Psi)(y,\th)\big|\,:\,\, (y,\th)\in\R\times
(0,1)\right\}<\infty.
\end{equation}
\item We denote by $\big\{\eps_{j,k}\,:\,\,(j,k)\in\Z^2\big\}$ the sequence of the real-valued independent ${\cal N}(0,1)$ Gaussian random variables defined,
for all $(j,k)\in\Z^2$, as,
\begin{equation}
\label{eq1:epsilonjk} \eps_{j,k}=2^{j/2} \int_{\R}\psi(2^j x-k)\,dW(x).
\end{equation}
\end{itemize}

Roughly speaking, the field $\{\widetilde{B}(t,\th)\}_{(t,\th)\in\R\times (0,1)}$ is defined as, 
$$
\widetilde{B}(t,\th)=\sum_{(j,k)\in\Z^2} 2^{-j\th}\eps_{j,k}\big [\Psi (2^j t-k,\th)-\Psi (-k,\th)\big].
$$
In order to precisely explain its definition, one needs some preliminary results. The following lemma allows to almost surely control the increase of the 
sequence $\big\{\eps_{j,k}\,:\,\,(j,k)\in\Z^2\big\}$.

\begin{lemma}
\label{lem:epsilonjk}
There is $\Omega^*$ an event of probability $1$, such
that every $\o\in\Omega^*$ satisfies the following two properties.
\begin{itemize}
\item[(i)] There exists $C$ a positive random variable, non depending on $(j,k)$ and of finite moment of any order, such that for all $(j,k)\in \Z^2$, one has
\begin{equation}
\label{eq2:epsilonjk}
|\eps_{j,k}(\o)|\le C(\o)\sqrt{\log(2+|j|)\log(2+|k|)}.
\end{equation}
\item[(ii)] For each fixed $s\in\R$ and $j\in\N$, let $\tau_j(s)$ be the random variable defined as,
\begin{equation}
\label{eq3:epsilonjkbis}
\tau_j (s)=\max\big\{|\eps_{j,k}|\,:\,\,k\in\Z \mbox{ and } |s-2^{-j}k|\le j 2^{1-j}\big\};
\end{equation}
then one has,
\begin{equation}
\label{eq3:epsilonjk}
\liminf_{j\rightarrow+\infty} \tau_j (s,\o)\ge 1/4.
\end{equation}
\end{itemize}
\end{lemma}

The proof of Lemma~\ref{lem:epsilonjk} has been omitted, since Part $(i)$ can be obtained similarly to Lemma 4 in \cite{ayache02} and Part $(ii)$ similarly to Lemma 4.1 in \cite{AJT07}. The precise definition of the field $\{\widetilde{B}(t,\th)\}_{(t,\th)\in\R\times (0,1)}$ is provided by the following proposition.

\begin{proposition}
\label{prop:absconc}
Let $\Psi$ be the function introduced in (\ref{eq1:pfpsi}) and let $\big\{\eps_{j,k}\,:\,\,(j,k)\in\Z^2\big\}$ be the sequence of the real-valued independent ${\cal N}(0,1)$ Gaussian random variables defined in (\ref{eq1:epsilonjk}). For all fixed $\o\in\Omega^*$ and $(t,\th)\in\R\times (0,1)$, one has 
\begin{equation}
\label{eq1:absconc}
\sum_{(j,k)\in\Z^2} 2^{-j\th}\big|\eps_{j,k}(\o)\big|\big|\Psi (2^j t-k,\th)-\Psi (-k,\th)\big|<\infty.
\end{equation}
Therefore, the series of real numbers 
$$
\sum_{(j,k)\in\Z^2} 2^{-j\th}\eps_{j,k}(\o)\big [\Psi (2^j t-k,\th)-\Psi (-k,\th)\big],
$$
converges to a finite limit which does not depend on the way the terms of the series are ordered; this limit is denoted by $\widetilde{B}(t,\th,\o)$.
Moreover for each $\o\notin \Omega^*$ and every $(t,\th)\in\R\times (0,1)$, one sets  $\widetilde{B}(t,\th,\o)=0$.
\end{proposition}

\noindent {\sc Proof of Proposition~\ref{prop:absconc}:} Let $\o\in\Omega^*$ and $(t,\th)\in\R\times (0,1)$ be arbitrary and fixed.
By using the triangle inequality, (\ref{eq1:prop-fdpsi}) in which one takes $m=n=0$, and (\ref{eq2:epsilonjk}), it follows that for all arbitrary fixed $j\in\N$,
\begin{eqnarray}
\label{eq2:absconc}
&& \sum_{k\in\Z} \big|\eps_{j,k}(\o)\big|\big|\Psi (2^j t-k,\th)-\Psi (-k,\th)\big|\nonumber\\
&& \le C_1 (\o) \sqrt{\log(2+j)}\sum_{k\in\Z} \left(\frac{\sqrt{\log(2+|k|)}}{\big (2+|2^j t-k|\big)^{2}}+\frac{\sqrt{\log(2+|k|)}}{\big (2+|k|\big)^{2}}\right)\nonumber\\
&& \le C_2 (\o) \sqrt{\log(2+j)\log(2+2^j |t|)}\sum_{k\in\Z} \left(\frac{\sqrt{\log(2+|k|)}}{\big (2+|2^j t-[2^j t]-k|\big)^{2}}+\frac{\sqrt{\log(2+|k|)}}{\big (2+|k|\big)^{2}}\right),\nonumber\\
\end{eqnarray}
where $[2^j t]$ denotes the integer part of $2^j t$ and where $C_1(\o)$ and $C_2(\o)$ are two finite constants non depending on $k$, $j$ and $t$. Then, noticing that,
\begin{equation}
\label{eq3:absconc}
\sup_{y\in [0,1]}\left\{ \sum_{k\in\Z}\frac{\sqrt{\log(2+|k|)}}{\big (2+|y-k|\big)^{2}}\right\}<\infty,
\end{equation}
it follows from (\ref{eq2:absconc}) that 
\begin{equation}
\label{eq4:absconc}
\sum_{j\in\N} \sum_{k\in\Z} 2^{-j\th}\big|\eps_{j,k}(\o)\big|\big|\Psi (2^j t-k,\th)-\Psi (-k,\th)\big|<\infty.
\end{equation}
Let us now prove that,
\begin{equation}
\label{eq5:absconc}
\sum_{j\in\Z_{-}} \sum_{k\in\Z} 2^{-j\th}\big|\eps_{j,k}(\o)\big|\big|\Psi (2^j t-k,\th)-\Psi (-k,\th)\big|<\infty.
\end{equation}
Assume that $j\in\Z_{-}$ is arbitrary and fixed. Applying, for any fixed $k\in\Z$, the Mean Value Theorem to the function $y\mapsto\Psi(y-k)$ on the interval $\big[\min\{2^{j} t, 0\}, \max\{2^{j} t, 0\}\big]\subseteq \big [-|t|,|t|\big ]$, one has that there exists $u\in \big [-|t|,|t|\big ]$, such that,
\begin{equation}
\label{eq6:absconc}
\Psi (2^j t-k,\th)-\Psi (-k,\th)=2^{j} t (\partial_{y} \Psi)(u-k,\th).
\end{equation}
Next, denote by ${\cal K}(t)$ and ${\cal K}^c(t)$ the sets ${\cal K}(t)=\{k\in\Z : |k|\le |t|\}$ and ${\cal K}^c(t)=\Z\setminus {\cal K}(t)$; observe that the cardinality of ${\cal K}(t)$ is bounded from above by $2|t|+1$. Putting together (\ref{eq6:absconc}), (\ref{eq1:prop-fdpsi}) in which one takes $(m,n)=(1,0)$, and (\ref{eq2:epsilonjk}), one obtains that,
\begin{eqnarray}
\label{eq7:absconc}
&& \sum_{k\in {\cal K}(t)} \big|\eps_{j,k}(\o)\big|\big|\Psi (2^j t-k,\th)-\Psi (-k,\th)\big|\nonumber\\
&& \le C_3(\o)\left(\sup_{y\in\R} \big|(\partial_{y} \Psi)(y,\th)\big| \right) 2^{j}|t|\sqrt{\log(2+|j|)} \sum_{k\in {\cal K}(t)}\sqrt{\log(2+|k|)} \nonumber\\
&& \le C_3(\o)|t|\big (2 |t|+1\big)\sqrt{\log(2+|t|)}\left(\sup_{y\in\R} \big|(\partial_{y} \Psi)(y,\th)\big|
\right) 2^{j} \sqrt{\log(2+|j|)}
\end{eqnarray}
and 
\begin{eqnarray}
\label{eq8:absconc}
&& \sum_{k\in {\cal K}^c(t)} \big|\eps_{j,k}(\o)\big|\big|\Psi (2^j t-k,\th)-\Psi (-k,\th)\big|\nonumber\\
&& \le C_3 (\o) 2^{j}|t|\sqrt{\log(2+|j|)}\sum_{k\in {\cal K}^c(t)}\left(\sup_{y\in [-|t|,|t|]} \big|(\partial_{y} \Psi)(y-k,\th)\big| \right) \sqrt{\log(2+|k|)} \nonumber\\
&& \le C_4 (\o) 2^{j}|t|\sqrt{\log(2+|j|)}\sum_{k=[|t|]+1}^{+\infty}\frac{\sqrt{\log(2+k)}}{\big(1+k-[|t|]\big)^{2}}\nonumber\\ 
&& \le C_5(\o) |t|\sqrt{\log(2+|t|)}\left(\sum_{k=0}^{+\infty}\frac{\sqrt{\log(2+k)}}{\big(2+k\big)^{2}}\right) 2^{j}\sqrt{\log(2+|j|)},
\end{eqnarray}
where $C_3 (\o)$ is the positive finite constant $C(\o)$ in (\ref{eq2:epsilonjk}), and where $C_4 (\o)$ and $C_5 (\o)$ are two positive finite constants non depending on $k$, $j$ and $t$. Next combining (\ref{eq7:absconc}) and (\ref{eq8:absconc}), with
the fact that $\th\in (0,1)$, one gets (\ref{eq5:absconc}). Finally (\ref{eq4:absconc}) and (\ref{eq5:absconc}) show that (\ref{eq1:absconc}) holds.
$\Box$
\\

Let us now explain the reason why the random field $\{\widetilde{B}(t,\th)\}_{(t,\th)\in\R\times (0,1)}$ can be identified with the random field 
$\{B(t,\th)\}_{(t,\th)\in\R\times (0,1)}$ defined in (\ref{eq:fieldB}).

\begin{proposition}
\label{prop:repwav-B}
The field $\{\widetilde{B}(t,\th)\}_{(t,\th)\in\R\times (0,1)}$ introduced in Proposition~\ref{prop:absconc}, is a modification of the field 
$\{B(t,\th)\}_{(t,\th)\in\R\times (0,1)}$ defined in (\ref{eq:fieldB}). 
\end{proposition}

\noindent {\sc Proof of Proposition~\ref{prop:repwav-B}:} The proof is quite classical in the area of multifractional processes (see e.g. \cite{BJR97,AT03,AJT07,ASX}). First one expands for all fixed $(t,\th)\in\R\times (0,1)$ the function $x\mapsto(t-x)_{+}^{\th-1/2}-(-x)_{+}^{\th-1/2}$ in the basis ${\cal ML}$, then one makes, in the deterministic integrals corresponding to the coefficients, the change of variable $u=2^j x-k$, finally using the isometry property of the Wiener integral in (\ref{eq:fieldB}), one can show that the series 
$$
\sum_{(j,k)\in\Z^2} 2^{-j\th}\eps_{j,k}\big [\Psi (2^j t-k,\th)-\Psi (-k,\th)\big],
$$
converges in $L^2(\Omega)$ ($\Omega$ is the underlying probability space) to the random variable $B(t,\th)$.
$\Box$\\

Let us now notice that for all $\o\in\Omega^*$ and $(t,\th)\in\R\times (0,1)$, $\widetilde{B}(t,\th,\o)$ can be expressed as,
\begin{equation}
\label{eq:decompBtilde}
\widetilde{B}(t,\th,\o)=\widetilde{B}_1(t,\th,\o)+\widetilde{B}_2(t,\th,\o)-\widetilde{R}(\th,\o),
\end{equation}
where 
\begin{equation}
\label{eq:Btilde1}
\widetilde{B}_1(t,\th,\o)=\sum_{j=-\infty}^{0} 2^{-j\th}\sum_{k\in\Z} \eps_{j,k}(\o)\big [\Psi (2^j t-k,\th)-\Psi (-k,\th)\big],
\end{equation}
\begin{equation}
\label{eq:Btilde2}
\widetilde{B}_2(t,\th,\o)=\sum_{j=1}^{+\infty} 2^{-j\th} \sum_{k\in\Z} \eps_{j,k}(\o)\Psi (2^j t-k,\th),
\end{equation}
and
\begin{equation}
\label{eq:procR}
\widetilde{R}(\th,\o)=\sum_{j=1}^{+\infty} 2^{-j\th}\sum_{k\in\Z} \eps_{j,k}(\o)\Psi (-k,\th).
\end{equation}
In the sequel, we show that the functions $(t,\th)\mapsto \widetilde{B}_1(t,\th,\o)$ and  $(t,\th)\mapsto \widetilde{R}(\th,\o)$ are $C^\infty$ over 
$\R\times (0,1)$; thus it turns out that for proving Theorem~\ref{theo:fieldB}, it is sufficient to show that it is true when the field $\{\widetilde{B}(t,\th)\}
_{(t,\th)\in\R\times (0,1)}$, is replaced by the more simple field $\{\widetilde{B}_2(t,\th)\}_{(t,\th)\in\R\times (0,1)}$. The following technical lemma will play a 
crucial r\^ole in the sequel.

\begin{lemma}
\label{lem:foncS}
For all $\o\in\Omega^*$ and $j\in\Z$, we denote by $S_j(\cdot,\cdot,\o)$ the real-valued function defined for every $(t,\th)\in\R\times (0,1)$ , as,
\begin{equation}
\label{eq1:foncS}
S_j(y,\th,\o)=\sum_{k\in\Z} \eps_{j,k}(\o)\Psi(y-k,\th).
\end{equation}
Then, the following two results hold, for each $\o\in\Omega^*$.
\begin{itemize}
\item[(i)] For all $j\in\Z$, the function $S_j(\cdot,\cdot,\o)$ is $C^\infty$ over $\R\times (0,1)$ and one has for every nonnegative integers $m,n$ and each $(y,\th)\in \R\times (0,1)$,
\begin{equation}
\label{eq1:foncSbis}
(\partial_{y}^{m} \partial_{\th}^{n}S_j)(y,\th,\o)=\sum_{k\in\Z} \eps_{j,k}(\o)(\partial_{y}^{m} \partial_{\th}^{n}\Psi)(y-k,\th).
\end{equation}
\item[(ii)] For all real numbers $0<a<b<1$ and for every nonnegative integers $m,n$, there exists a finite constant $C(\o)$ only depending on 
$a,b,m,n,\o$, such that for all $j\in\Z$ and positive real number $z$, one has,
\begin{eqnarray}
\label{eq2:foncS}
&& \sup\left\{ \big|(\partial_{y}^{m} \partial_{\th}^{n}S_j)(y,\th,\o)\big| \,:\,\, (y,\th)\in [-z,z]\times [a,b]\right\}\nonumber\\
&& \le  \sup\left\{\sum_{k\in\Z} |\eps_{j,k}(\o)|\big|(\partial_{y}^{m} \partial_{\th}^{n}\Psi\big)(y-k,\th)\big| \,:\,\, (y,\th)\in [-z,z]\times [a,b]\right\}\nonumber\\
&& \le C(\o) \sqrt{\log(2+|z|)\log(2+|j|)}.
\end{eqnarray}
\end{itemize}
\end{lemma}

\noindent {\sc Proof of Lemma~\ref{lem:foncS}:} First observe that by using (\ref{eq2:epsilonjk}) and (\ref{eq1:prop-fdpsi}), one can prove that for each fixed
$(j,y,\th,\o)\in \Z\times\R\times (0,1)\times\Omega^*$, the series in (\ref{eq1:foncS}) is absolutely convergent, which implies that the function $S_j(\cdot,\cdot,\o)$
is well-defined. Let us now show that Part $(i)$ holds, to this end, it is sufficient to prove that for all nonnegative integers $m,n$, the series 
$$
\sum_{k\in\Z} \eps_{j,k}(\o)(\partial_{y}^{m} \partial_{\th}^{n}\Psi)(y-k,\th),
$$
is uniformly convergent, on each compact set of the form $[-M,M]\times [a,b]$, where $M,a,b$ are arbitrary real numbers satisfying $M>0$ and $0<a<b<1$. In order 
to derive the latter result, we will show that,
\begin{equation} 
\label{eq3:foncS}
\sum_{k\in\Z} |\eps_{j,k}(\o)| s_k <\infty,
\end{equation}
where for each $k\in\Z$,
$$
s_k=\sup\left\{ \big|(\partial_{y}^{m} \partial_{\th}^{n}\Psi)(y-k,\th)\big|\,:\,\, (y,\th)\in [-M,M]\times [a,b]\right\}.
$$
Observe that, in view of (\ref{eq1:prop-fdpsi}), one has, for some constant $c_1>0$ and all $k\in\Z$, satisfying $|k|>M$,
$$
s_k\le c_1 \big (2+|k|-M\big)^{-2}.
$$
Therefore using (\ref{eq2:epsilonjk}), one gets (\ref{eq3:foncS}). Let us now prove that Part $(ii)$ holds. It follows from (\ref{eq1:foncSbis}), (\ref{eq1:prop-fdpsi}) and (\ref{eq2:epsilonjk}), that, for all $j\in\Z$ and $(y,\th)\in [-z,z]\times [a,b]$,  
\begin{eqnarray}
\label{eq4:foncS}
&& \big|(\partial_{y}^{m} \partial_{\th}^{n}S_j)(y,\th,\o)\big|\nonumber\\
&& \le \sum_{k\in\Z} |\eps_{j,k}(\o)|\big|(\partial_{y}^{m} \partial_{\th}^{n}\Psi)(y-k,\th)\big|\nonumber\\
&& \le C_2 (\o) \sqrt{\log(2+|j|)} \sum_{k\in\Z} \frac{\sqrt{\log(2+|k|)}}{\big(2+|y-k|\big)^{2}}\nonumber\\
&& \le C_3 (\o) \sqrt{\log(2+|j|)\log(2+z)}\sum_{k\in\Z} \frac{\sqrt{\log(2+|k|)}}{\big(2+|y-[y]-k|\big)^{2}},
\end{eqnarray}
where $[y]$ is the integer part of $y$ and where $C_2(\o)$ and $C_3(\o)$ are two constants non depending on $j$, $y$, $z$ and $\th$. Then, combining (\ref{eq4:foncS})
with (\ref{eq3:absconc}), one obtains (\ref{eq2:foncS}).
$\Box$\\

The following lemma corresponds to Proposition~2.1 in \cite{AT03}, yet we prefer to give a short proof of it, for the sake of clarity.

\begin{lemma}
\label{lem:infClf}
For all $\o\in\Omega^*$ the function $(t,\th)\mapsto \widetilde{B}_1(t,\th,\o)$ (see (\ref{eq:Btilde1})) is $C^\infty$ over $\R\times (0,1)$.
\end{lemma}

\noindent {\sc Proof of Lemma~\ref{lem:infClf}:} In view of (\ref{eq:Btilde1}) and (\ref{eq1:foncS}), one has that 
$$
\widetilde{B}_1(t,\th,\o)=\sum_{j=-\infty}^{0} Q_j(t,\th,\o),
$$
where for each $j\in\Z_{-}$, $Q_j(\cdot,\cdot,\o)$ is the $C^\infty$ function over $\R\times (0,1)$, defined as,
\begin{equation}
\label{eq1:infClf}
Q_j(t,\th,\o)=2^{-j\th} \big[ S_j(2^j t,\th,\o)-S_j(0,\th,\o)\big].
\end{equation}
Therefore, in order to show the lemma, it is sufficient to prove that for each nonnegative integers $p,q$, and real numbers $M,a,b$ satisfying $M>0$ and 
$0<a<b<1$, one has,
\begin{equation}
\label{eq2:infClf}
\sum_{j=-\infty}^{0}\sup\left\{\big|(\partial_{t}^p \partial_{\th}^q Q_j)(t,\th,\o)\big|\,:\,\, (t,\th)\in [-M,M]\times [a,b]\right\}<\infty.
\end{equation}
Let us first study the case where $p=0$ and $q$ is an arbitrary nonnegative integer. Using (\ref{eq1:infClf}) and the Leibniz formula, it follows that
\begin{equation}
\label{eq3:infClf}
(\partial_{\th}^q Q_j)(t,\th,\o)=\sum_{l=0}^q \binom{q}{l} (-j\log 2)^l 2^{-j\th}\big[(\partial_{\th}^{q-l} S_j)(2^j t,\th,\o)-(\partial_{\th}^{q-l} S_j)(0,\th,\o)\big],
\end{equation}
where $\binom{q}{l}$ is the binomial coefficient $\frac{q!}{l!\times (q-l)!}$. Next applying the Mean Value Theorem, one gets that 
\begin{eqnarray}
\label{eq4:infClf}
&&\sup\left\{\big|(\partial_{\th}^{q-l} S_j)(2^j t,\th,\o)-(\partial_{\th}^{q-l} S_j)(0,\th,\o)\big|\,:\,\, (t,\th)\in [-M,M]\times [a,b]\right\}\nonumber\\
&& \le 2^{j}M\sup\left\{\big|(\partial_y \partial_{\th}^{q-l} S_j)(y,\th)\big|\,:\,\, (y,\th)\in [-M,M]\times [a,b]\right\}.
\end{eqnarray}
Next putting together (\ref{eq3:infClf}), (\ref{eq4:infClf}) and (\ref{eq2:foncS}), it follows that (\ref{eq2:infClf}) holds in the case where $p=0$. Let us now
study the case where $p\ge 1$ and $q$ is an arbitrary nonnegative integer. In view of (\ref{eq1:infClf}), one has that,
$$
(\partial_{t}^p Q_j)(t,\th,\o)=2^{j(p-\th)}(\partial_{y}^p S_j)(2^j t,\th,\o).
$$
Therefore, using the Leibniz formula, one obtains that,
\begin{equation}
\label{eq5:infClf}
(\partial_{t}^p \partial_{\th}^q Q_j)(t,\th,\o)=\sum_{l=0}^q \binom{q}{l} (-j\log 2)^l 2^{j(p-\th)}(\partial_{y}^p\partial_{\th}^{q-l} S_j)(2^j t,\th,\o).
\end{equation}
Finally, combining (\ref{eq5:infClf}) with (\ref{eq2:foncS}), one gets (\ref{eq2:infClf}). $\Box$
\\

The proof of the following lemma has been omitted since it is rather similar to that of Lemma~\ref{lem:infClf}.
\begin{lemma}
\label{lem:infCR}
For all $\o\in\Omega^*$ the function $(t,\th)\mapsto \widetilde{R}(\th,\o)$ (see (\ref{eq:procR})) is $C^\infty$ over $\R\times (0,1)$.
\end{lemma}

Let us now give some results concerning the global regularity of the function $(t,\th)\mapsto \widetilde{B}_2(t,\th,\o)$ (see (\ref{eq:Btilde2})).

\begin{lemma}
\label{lem:supinfClf}
For all $\o\in\Omega^*$, the following three results hold.
\begin{itemize}
\item[(i)] The function $(t,\th)\mapsto \widetilde{B}_2(t,\th,\o)$ is continuous over $\R\times (0,1)$.
\item[(ii)] For each fixed $t\in\R$, the function $\th\mapsto \widetilde{B}_2(t,\th,\o)$ is $C^\infty$ over $(0,1)$.
\item[(iii)] For every fixed nonnegative integer $n$, the function $(t,\th)\mapsto (\partial_{\th}^n \widetilde{B}_2)(t,\th,\o)$ is continuous over $\R\times (0,1)$. 
\end{itemize}
\end{lemma}

\noindent {\sc Proof of Lemma~\ref{lem:supinfClf}:} In view of (\ref{eq:Btilde2}) and (\ref{eq1:foncS}), one has that 
\begin{equation}
\label{eq1:supinfClf}
\widetilde{B}_2(t,\th,\o)=\sum_{j=1}^{+\infty} V_j(t,\th,\o),
\end{equation}
where for each $j\in\N$, $V_j(\cdot,\cdot,\o)$ is the $C^\infty$ function over $\R\times (0,1)$, defined as,
\begin{equation}
\label{eq2:supinfClf}
V_j(t,\th,\o)=2^{-j\th}S_j(2^j t,\th,\o).
\end{equation}
Therefore, in order to show the lemma, it is sufficient to prove that for each nonnegative integer $q$, and real numbers $M,a,b$ satisfying $M>0$ and 
$0<a<b<1$, one has,
\begin{equation}
\label{eq3:supinfClf}
\sum_{j=1}^{+\infty}\sup\left\{\big|(\partial_{\th}^q V_j)(t,\th,\o)\big|\,:\,\, (t,\th)\in [-M,M]\times [a,b]\right\}<\infty.
\end{equation}
By using (\ref{eq2:supinfClf}) and the Leibniz formula, it follows that,
\begin{equation}
\label{eq4:supinfClf}
(\partial_{\th}^q V_j)(t,\th,\o)=\sum_{l=0}^q \binom{q}{l} (-j\log 2)^l 2^{-j\th}(\partial_{\th}^{q-l} S_j)(2^j t,\th,\o).
\end{equation}
Finally, combining (\ref{eq4:supinfClf}) with (\ref{eq2:foncS}), one gets (\ref{eq3:supinfClf}).
$\Box$\\

Now we are in position to prove Parts $(i)$ and $(iii)$ of Theorem~\ref{theo:fieldB}.

\noindent {\sc Proof of Parts $(i)$ and $(iii)$ of Theorem~\ref{theo:fieldB}:} These two parts are straightforward consequences of (\ref{eq:decompBtilde}), and Lemmas~\ref{lem:infClf}, \ref{lem:infCR} and \ref{lem:supinfClf}. $\Box$\\

Let us now turn to the proof of Theorem~\ref{theo:fieldB} Part $(iv)$. We need the following two lemmas.

\begin{lemma} 
\label{lem:Lpschunif}
For all $\o\in\Omega^*$, for each nonnegative integer $n$, and for every real numbers $M,a,b$ satisfying $M>0$ and $0<a<b<1$, one 
has for all $\th_1,\th_2\in [a,b]$,
\begin{equation}
\label{eq:Lpschunif}
\sup_{t\in [-M,M]} \big| (\partial_{\th}^n \widetilde{B}_2)(t,\th_1,\o)-(\partial_{\th}^n \widetilde{B}_2)(t,\th_2,\o)\big| \le C(\o) |\th_1-\th_2|,
\end{equation}
where the finite constant $C(\o)$, only depends on $\o,n,M,a,b$.
\end{lemma}

Observe that in the case where $n=0$, Lemma~\ref{lem:Lpschunif} can be related to Part $(c)$ of Proposition~2.2 in \cite{AT03}.\\

\noindent {\sc Proof of Lemma~\ref{lem:Lpschunif}:} The lemma easily results from the Mean Value Theorem and Lemma~\ref{lem:supinfClf}. $\Box$

\begin{lemma}
\label{lem:part-modcB}
For all $\o\in\Omega^*$, for each nonnegative integer $n$, for every arbitrarily small real number $\eps>0$ and for every real numbers $M,a,b$ satisfying $M>0$ and $0<a<b<1$, one has, for each  $\th\in [a,b]$ and $t_1,t_2\in [-M,M]$, 
\begin{equation}
\label{eq1:part-modcB} 
\big |(\partial_{\th}^n \widetilde{B}_2)(t_1,\th,\o)-(\partial_{\th}^n \widetilde{B}_2)(t_2,\th,\o)\big |\le C_1(\o) |t_1-t_2|^{\th-\eps}.
\end{equation}
where the finite constant $C(\o)$, only depends on $\o,n,\eps,M,a,b$.
\end{lemma}
\noindent {\sc Proof of Lemma~\ref{lem:part-modcB}:} The proof is inspired by that of Proposition~4.2 in \cite{AJT07}. First observe that, in view of (\ref{eq1:supinfClf}) and (\ref{eq4:supinfClf}), in order to derive 
(\ref{eq1:part-modcB}), it is sufficient to show that for all $l\in\{0,\ldots, n\}$, one has,
\begin{equation}
\label{eq5:part-modcB} 
\sum_{j=1}^{+\infty} j^l 2^{-j\th}\big|(\partial_{\th}^{n-l} S_j)(2^j t_1,\th,\o)-(\partial_{\th}^{n-l} S_j)(2^j t_2,\th,\o)\big| \le C_1(\o) |t_1-t_2|^{\th-\eps},
\end{equation}
where $C_1(\o)$ is a finite constant only depending on $\o,n,\eps,M,a,b$. Also observe that there exists a constant $c_2>0$, only depending, on $\eps$, $M$ and $n$ such that,
for all integer $j\ge 1$,
\begin{equation}
\label{eq3:part-modcB} 
j^l 2^{-j\th}\sqrt{\log(2+M2^j)\log(2+|j|)}\le c_2 2^{-j(\th-\eps)}.
\end{equation}
The inequality (\ref{eq5:part-modcB}) is clearly satisfied when $t_1=t_2$, so from now on we assume that $|t_1-t_2|>0$. Let then $j_0\ge 1$ be the unique integer such that
\begin{equation}
\label{eq2:part-modcB}
M 2^{-j_0 +1}<|t_1-t_2|\le M 2^{-j_0 +2}.
\end{equation}
By using (\ref{eq2:foncS}), (\ref{eq3:part-modcB}) and (\ref{eq2:part-modcB}), one has that,
\begin{eqnarray}
\label{eq4:part-modcB}
&&\sum_{j=j_0+1}^{+\infty} j^l 2^{-j\th}\left(\sup\left\{\big|(\partial_{\th}^{n-l} S_j)(2^j t,\th,\o)\big|\,:\,\,(t,\th)\in [-M,M]\times [a,b]\right\}\right)\\
&& \le C_3(\o)\sum_{j=j_0+1}^{+\infty} 2^{-j(\th-\eps)}\le C_3(\o)  \big (1-2^{-(a-\eps)}\big)^{-1} 2^{-(j_0+1)(\th-\eps)}\le C_4(\o)|t_1-t_2|^{\th-\eps},\nonumber
\end{eqnarray}
where $C_3(\o)$ and $C_4(\o)$ are two finite constants only depending on $\o,n,\eps,M,a,b$. On the other hand, 
it follows from the Mean Value Theorem, that, for all $j\in\N$,
\begin{eqnarray}
\label{eq4bis:part-modcB}
&& \big|(\partial_{\th}^{n-l} S_j)(2^j t_1,\th,\o)-(\partial_{\th}^{n-l} S_j)(2^j t_2,\th,\o)\big|\\
&& \le 2^{j}|t_1-t_2| \left(\sup\left\{\big|(\partial_y \partial_{\th}^{n-l} S_j)(y,\th,\o)\big|\,:\,\,(y,\th)\in [-2^{j}M,2^{j} M]\times [a,b]\right\}\right)\nonumber.
\end{eqnarray}
Then, putting together, (\ref{eq4bis:part-modcB}), (\ref{eq2:foncS}), (\ref{eq3:part-modcB}) and (\ref{eq2:part-modcB}), one obtains that,
\begin{eqnarray}
\label{eq6:part-modcB}
&&\sum_{j=1}^{j_0} j^l 2^{-j\th}\big|(\partial_{\th}^{n-l} S_j)(2^j t_1,\th,\o)-(\partial_{\th}^{n-l} S_j)(2^j t_2,\th,\o)\big|\nonumber\\
&& \le C_5 (\o) |t_1-t_2|\sum_{j=1}^{j_0} 2^{(1-\th+\eps)j}\le C_5 (\o) |t_1-t_2| \big(2^{(1-b+\eps)}-1\big)^{-1} 2^{(1-\th+\eps)(j_0+1)}\nonumber\\
&& \le C_6 (\o) |t_1-t_2|^{\th-\eps},
\end{eqnarray}
where $C_5(\o)$ and $C_6(\o)$ are two finite constants only depending on $\o,n,\eps,M,a,b$. Finally combining (\ref{eq4:part-modcB}) with (\ref{eq6:part-modcB}), one gets (\ref{eq5:part-modcB}).
$\Box$
\\

Now we are in position to prove Part $(iv)$ of Theorem~\ref{theo:fieldB}.

\noindent {\sc Proof of Part $(iv)$ of Theorem~\ref{theo:fieldB}:} it follows from (\ref{eq:decompBtilde}),  Lemma~\ref{lem:infClf} and Lemma~\ref{lem:infCR}, that 
it is sufficient to show that (\ref{eq3:TfieldB}) holds, when $\partial_{\th}^n \widetilde{B}$ is replaced by $\partial_{\th}^n \widetilde{B}_2$. Let 
$(t_1,\th_1)\in [-M,M]\times [a,b]$ and $(t_2,\th_2)\in [-M,M]\times [a,b]$, there is 
no restriction to assume that $\th_1=\max\{\th_1,\th_2\}$. Using the triangle inequality, one has that, 
\begin{eqnarray*}
&& \big|(\partial_{\th}^n \widetilde{B}_2)(t_1,\th_1,\o)-(\partial_{\th}^n \widetilde{B}_2)(t_2,\th_2,\o)\big |\\
&& \le \big|(\partial_{\th}^n \widetilde{B}_2)(t_1,\th_1,\o)-(\partial_{\th}^n \widetilde{B}_2)(t_2,\th_1,\o)\big |+\big|(\partial_{\th}^n \widetilde{B}_2)(t_2,\th_1,\o)-(\partial_{\th}^n \widetilde{B}_2)(t_2,\th_2,\o)\big|.
\end{eqnarray*}
Next, combining the latter inequality with Lemmas~\ref{lem:Lpschunif}~and~\ref{lem:part-modcB}, we can finish our proof. $\Box$\\

Let us now turn to the proof of Part $(ii)$ of Theorem~\ref{theo:fieldB}. Notice that (\ref{eq1:TfieldB}) is a straightforward consequence of (\ref{eq3:TfieldB}) in
which one takes $n=0$ and $M,a,b$ such that $(s,\th)\in (-M,M)\times [a,b]$. In order to show that (\ref{eq2:TfieldB}) holds, one needs to introduce the real-valued function $\widetilde{\Psi}$, defined for every
$(y,\th)\in\R\times (0,1)$ as,
\begin{equation}
\label{eq:dfpsi}
\widetilde{\Psi}(y,\th)=\frac{1}{2\pi\Ga(\th+1/2)}\int_{\R}e^{iy\xi}e^{-i\,\sgn(\xi) (\th+1/2)\frac{\pi}{2}}|\xi|^{\th+1/2}\widehat{\psi}(\xi)\,d\xi,
\end{equation}
where $\widehat{\psi}$ is the Fourier transform of the Lemari\'e-Meyer mother wavelet $\psi$ introduced at the very beginning of this subsection. Let us now 
give some useful properties of $\widetilde{\Psi}$. The proof of the following lemma, has been omitted since it is very similar to those of Lemmas 2.1 and 3.4 in \cite{AT03}.

\begin{lemma}
\label{lem:prop-fdpsi}
\begin{itemize}
\item[(i)] $\widetilde{\Psi}$ is $C^\infty$ over $\R\times (0,1)$ and its partial derivatives
of any order are well-localized in the variable $y\in\R$, uniformly
in the variable $\th\in (0,1)$; in other words,  one has, for every nonnegative integers $m$ and $n$,
\begin{equation}
\label{eq1:prop-fdpsitilde} \sup\left\{\big(2+|y|\big)^2
\big|(\partial_{y}^m \partial_{\th}^n
\widetilde{\Psi})(y,\th)\big|\,:\,\, (y,\th)\in\R\times
(0,1)\right\}<\infty.
\end{equation}
\item[(ii)] For all $\th\in (0,1)$, the first moment of the function $\widetilde{\Psi}(\cdot,\th)$ vanishes, that is
\begin{equation}
\label{eq2:prop-fdpsi}
\int_\R \widetilde{\Psi}(y,\th)\,dy=0.
\end{equation}
\item[(iii)] For all $\th\in (0,1)$, the system of functions
$ \big\{2^{j'/2}\Psi (2^{j'}\cdot-k',\th)\,:\,\,j'\in\Z,\,k'\in\Z\big\}$ (recall that $\Psi$ has been introduced in (\ref{eq1:pfpsi}))
and $ \big\{2^{j/2}\widetilde{\Psi}
(2^j\cdot-k,\th)\,:\,\,j\in\Z,\,k\in\Z\big\} $ is biorthogonal. This
means that for all $j,k,j',k'\in\Z$, one has
\begin{equation}
\label{eq3:prop-fdpsi}
2^{(j+j')/2}\int_{\R}\Psi (2^{j'} t-k',\th)\widetilde{\Psi} (2^{j}t-k,\th)\,dt=\de (j,k;j',k'),
\end{equation}
where $\de (j,k;j',k')=1$ if $(j,k)=(j',k')$ and $0$ otherwise.
\end{itemize}
\end{lemma}
The following lemma, which can be related to Part $(c)$ of Proposition~3.3 in \cite{AJT07}, allows one to understand the motivation behind the introduction of $\widetilde{\Psi}$.
\begin{lemma}
\label{lem:wavtran}
For all $\o\in\Omega^*$, let $\widetilde{B}_2(\cdot,\cdot,\o)$ be the function introduced in (\ref{eq:Btilde2}). Then for each fixed $(j,k,\th, \o)\in 
\N\times\Z\times (0,1)\times \Omega^*$, the integral,
\begin{equation}
\label{eq1:wavtran}
\widetilde{{\cal W}}_{j,k}(\th,\o)=2^{j(1+\th)}\int_{\R} \widetilde{B}_2(s,\th,\o)\widetilde{\Psi}(2^j s-k,\th)\,ds,
\end{equation}
is well-defined. Moreover, one has,
\begin{equation}
\label{eq2:wavtran}
\widetilde{{\cal W}}_{j,k}(\th,\o)=\eps_{j,k}(\o),
\end{equation}
where $\eps_{j,k}$ is the ${\cal N}(0,1)$ Gaussian random variable introduced in (\ref{eq1:epsilonjk}).
\end{lemma}

\noindent {\sc Proof of Lemma~\ref{lem:wavtran}:} First observe that by using the second inequality in (\ref{eq2:foncS}), in the case where $m=n=0$, $z=1+|s|$ and $\th\in [a,b]$, one obtains,
in view of (\ref{eq:Btilde2}), that 
\begin{eqnarray}
\label{eq3:wavtran}
\big|\widetilde{B}_2(s,\th,\o)\big|&\le & \sum_{j'=1}^{+\infty} 2^{-j'\th} \sum_{k'\in\Z} |\eps_{j',k'}(\o)\| \big|\Psi (2^{j'} s-{k'},\th)\big|\nonumber\\
&\le & C_1 (\o)\sqrt{\log(3+|s|)},
\end{eqnarray}
where $C_1(\o)$ is a constant only depending on $a,b,\o$. Thus (\ref{eq3:wavtran}) and (\ref{eq1:prop-fdpsitilde}) imply that the integral in (\ref{eq1:wavtran}) is well-defined. Moreover (\ref{eq:Btilde2}), the dominated convergence Theorem and (\ref{eq3:prop-fdpsi}), entail that (\ref{eq2:wavtran}) holds. $\Box$\\
 
Now, we are in position to prove Part $(ii)$ of Theorem~\ref{theo:fieldB}.

\noindent {\sc Proof of Part $(ii)$ of Theorem~\ref{theo:fieldB}:} As we have mentioned before, Relation
(\ref{eq1:TfieldB}) is a straightforward
consequence of Part $(iv)$ of the theorem, which has already been proved. So, it remains to show that
(\ref{eq2:TfieldB}) holds; the proof, we are going to give, is inspired by that of Proposition~4.1 in \cite{AJT07}. It follows from (\ref{eq:decompBtilde}),  Lemma~\ref{lem:infClf} and Lemma~\ref{lem:infCR}, that 
it is sufficient to show that (\ref{eq2:TfieldB}) holds, when $\widetilde{B}$ is replaced by $\widetilde{B}_2$. Suppose {\it ad absurdum} that the latter
relation is not satisfied for some $\o_0\in\Omega^*$,
$(s_0,\th_0)\in\R\times (0,1)$ and $\eps_0>0$ (notice that there is no restriction to assume that $\th_0+\eps_0<1$), then there exists a finite
constant $C_1 (\o_0)>0$ such that for some finite constant $\eta_0>0$ and for all real number $s$ satisfying $|s-s_0|\le\eta_0$, one
has
\begin{equation}
\label{eq4:TfieldB}
\big |\widetilde{B}_2(s,\th_0,\o_0)-\widetilde{B}_2(s_0,\th_0,\o_0)\big|\le C_1 (\o_0)|s-s_0|^{\th_0+\eps_0}.
\end{equation}
On the other hand, by using the fact that $s\mapsto \big |\widetilde{B}_2(s,\th_0,\o_0)-\widetilde{B}_2(s_0,\th_0,\o_0)\big|$ is a continuous 
function over $\R$ as well as (\ref{eq3:wavtran}), one obtains that,
\begin{equation}
\label{eq5:TfieldB}
\sup\left\{\frac{\big |\widetilde{B}_2(s,\th_0,\o_0)-\widetilde{B}_2(s_0,\th_0,\o_0)\big|}{|s-s_0|^{\th_0+\eps_0}}\,:\,\, s\in\R\mbox{ and } |s-s_0|\ge\eta_0\right\}
< \infty.
\end{equation}
Thus, combining (\ref{eq4:TfieldB}) with (\ref{eq5:TfieldB}), it follows that, there is a finite constant $C_2(\o_0)>0$, such that for all $s\in\R$,
\begin{equation}
\label{eq6:TfieldB}
\big |\widetilde{B}_2(s,\th_0,\o_0)-\widetilde{B}_2(s_0,\th_0,\o_0)\big|\le C_2 (\o_0)|s-s_0|^{\th_0+\eps_0}.
\end{equation}
Then using (\ref{eq1:wavtran}), (\ref{eq2:wavtran}), (\ref{eq2:prop-fdpsi}), (\ref{eq6:TfieldB}) and the change of variable $t=2^j s-k$, one gets that for all integer $j\ge 1$ and $k\in\Z$,
\begin{eqnarray}
\label{eq7:TfieldB}
\nonumber
|\eps_{j,k}(\o_0)|&=& 2^{j(1+\th_0)}\Big |\int_{\R} \widetilde{B}_2(s,\th_0,\o_0)\widetilde{\Psi}(2^j s-k,\th_0)\,ds\Big |\\
\nonumber
&=& 2^{j(1+\th_0)}\Big |\int_{\R} \big [\widetilde{B}_2(s,\th_0,\o_0)-\widetilde{B}_2(s_0,\th_0,\o_0)\big ]\widetilde{\Psi}(2^j s-k,\th_0)\,ds\Big |\\
\nonumber
&\le & C_2(\o_0) 2^{j(1+\th_0)}\int_{\R} |s-s_0|^{\th_0+\eps_0}\big|\widetilde{\Psi}(2^j s-k,\th_0)\big|\,ds\\
\nonumber
& \le & C_2(\o_0) 2^{j\th_0}\int_{\R}\big |2^{-j}t+2^{-j}k-s_0\big|^{\th_0+\eps_0}|\widetilde{\Psi}(t,\th_0)|\,dt\\
&\le & C_3(\o_0) 2^{-j\eps_0}\Big (1+\big |2^j s_0 -k\big|^{\th_0+\eps_0}\Big ),
\end{eqnarray}
where $C_3(\o_0)$ is a finite constant non depending on $j$ and $k$.
Observe that to derive the last inequality in (\ref{eq7:TfieldB}), we have used the fact that, for some finite constant $c_4$ and for every real numbers
$u,v$, one has, $\big (|u|+|v|\big )^{\th_0+\eps_0}\le
c_4\big(|u|^{\th_0+\eps_0}+|v|^{\th_0+\eps_0}\big)$, and we have also used (\ref{eq1:prop-fdpsitilde}). Finally, (\ref{eq7:TfieldB}) and (\ref{eq3:epsilonjkbis}) entail that 
$$
\tau_j(s_0,\o_0)=\mathcal{O}(2^{-j\eps_0/2}),
$$
which contradicts (\ref{eq3:epsilonjk}). $\Box$

\subsection{Proof of Proposition~\ref{prop:rep-partialB}}
\label{subsec:rep-partialB}

Let us denote by $\{Z(t,\th)\}_{(t,\th)\in \R\times (0,1)}$ the field defined for each $(t,\th)$ as the Wiener integral,
\begin{equation}
\label{eq2:rep-partialB}
Z(t,\th)= \int_{\R}\Big\{(t-x)_{+}^{\th-1/2}\log\big [(t-x)_+\big ]-(-x)_{+}^{\th-1/2}\log\big [(-x)_+\big ]\Big\}\,dW(x).
\end{equation}
By expanding for every fixed $(t,\th)\in\R\times (0,1)$ the function
$$
x\mapsto (t-x)_{+}^{\th-1/2}\log\big [(t-x)_+\big ]-(-x)_{+}^{\th-1/2}\log\big [(-x)_+\big ],
$$
in the Lemari\'e-Meyer orthonormal wavelet basis ${\cal ML}$ (see the beginning of Subsection~\ref{subsec:wav}), and by using the isometry property of the Wiener integral, it follows that
\begin{equation}
\label{eq3:rep-partialB}
Z(t,\th)= \sum_{j\in\Z}\sum_{k\in\Z}\eps_{j,k}
\big [d_{j,k}(t,\th)-d_{j,k}(0,\th)\big ],
\end{equation}
where the ${\cal N}(0,1)$ Gaussian random variables $\eps_{j,k}$ have been introduced in (\ref{eq1:epsilonjk}) and where
\begin{equation}
\label{eq4:rep-partialB}
d_{j,k}(t,\th)=2^{j/2}\int_{\R} (t-x)_{+}^{\th-1/2}\log\big [(t-x)_+\big ]\psi(2^jx-k)\,dx.
\end{equation}
Observe that the series in (\ref{eq3:rep-partialB}) is convergent, for every fixed $(t,\th)$,
in $L^2(\Omega)$, where $\Omega$ is the underlying probability space. Setting in (\ref{eq4:rep-partialB}) $s=2^j x-k$, and using Lemma~\ref{lem:foncS}, one obtains after a sequence of standard computations, that
for every $\o\in\Omega^*$ (the event of probability $1$ introduced in Theorem~\ref{theo:fieldB}) and for each $(t,\th)\in\R\times (0,1)$, one has when $j\le 0$,
\begin{equation}
\label{eq5:rep-partialB}
(\partial_\th Q_j)(t,\th,\o)=\sum_{k\in\Z}\eps_{j,k}(\o)\big [d_{j,k}(t,\th)-d_{j,k}(0,\th)\big ].
\end{equation}
and, one has when $j>0$, 
\begin{equation}
\label{eq6:rep-partialB}
(\partial_\th V_j)(t,\th,\o)-(\partial_\th V_j)(0,\th,\o)=\sum_{k\in\Z}\eps_{j,k}(\o)\big [d_{j,k}(t,\th)-d_{j,k}(0,\th)\big ];
\end{equation}
recall that the $C^\infty$ functions $Q_j(\cdot,\cdot,\o)$ and $V_j(\cdot,\cdot,\o)$ have been introduced in (\ref{eq1:infClf}) and in (\ref{eq2:supinfClf}).
Moreover using Proposition~\ref{prop:absconc}, (\ref{eq2:infClf}) and (\ref{eq3:supinfClf}), it follows that for all $\o\in\Omega^*$ and $(t,\th)\in\R\times (0,1)$,
\begin{equation}
\label{eq7:rep-partialB}
(\partial_\th B)(t,\th,\o)=\sum_{j=-\infty}^{0} (\partial_\th Q_j)(t,\th,\o)+\sum_{j=1}^{+\infty}\big[(\partial_\th V_j)(t,\th,\o)-(\partial_\th V_j)(0,\th,\o)\big].
\end{equation}
Finally putting together (\ref{eq3:rep-partialB}), (\ref{eq5:rep-partialB}), (\ref{eq6:rep-partialB}) and (\ref{eq7:rep-partialB}), one obtains the proposition. $\Box$
\\

\noindent{\bf Acknowledgement. } The author thanks the anonymous referees and Professor Davar Khoshnevisan, for their
valuable comments, which have led to improvements of the manuscript. 
Some parts of it have been written, while the author was invited Professor at the department
of Mathematics of the Wuhan University in China; he is very grateful
to this department for its financial support, also he thanks all of its members for their kindness and in particular
Professor Yijun Hu.

\bibliographystyle{plain}
\begin{small}

\end{small}

\vskip.4in

\begin{quote}
\begin{small}

\textsc{Antoine Ayache}: U.M.R. CNRS 8524, Laboratoire Paul
Painlev\'e, B\^atiment M2, Universit\'e Lille 1, 59655 Villeneuve d'Ascq Cedex, France.\\
E-mail: \texttt{Antoine.Ayache@math.univ-lille1.fr}\\[1mm]

\end{small}
\end{quote}

\end{document}